\documentclass[11pt,letterpaper]{article}
\usepackage[T1]{fontenc}
\usepackage{lmodern}
\usepackage[margin=1in]{geometry}
\usepackage{amsmath,amssymb,amsthm}

\usepackage[colorlinks=true,linkcolor=blue]{hyperref}

\usepackage{graphicx}
\usepackage{float}
\usepackage[numbers,square,sort&compress]{natbib}
\usepackage{microtype}
\usepackage{booktabs,pgfplots}
\pgfplotsset{compat=1.18}
\allowdisplaybreaks[1]
\newtheorem{assumption}{Assumption}
\newtheorem{proposition}{Proposition}
\newtheorem{corollary}{Corollary}
\newtheorem{lemma}{Lemma}
\newtheorem{theorem}{Theorem}
\theoremstyle{remark}
\newtheorem{remark}{Remark}
\newcommand{\E}{\mathbb E}
\newcommand{\dd}{\,\mathrm d}
\newtheorem{definition}{Definition}

\title{Optimal Subsidy Design for Productive Investment across Markets}
\author{%
  \textsc{Alexandros Pavlis}\thanks{%
    Department of Statistics, London School of Economics and
    Political Science, Houghton Street, London WC2A 2AE, UK.
    \newline
    \texttt{a.pavlis@lse.ac.uk}%
  }%
}
\date{}

\begin{document}
\maketitle
\vspace{-2em}
\begin{abstract}
We study the optimal design of investment subsidies across heterogeneous monopolistic markets subject to a common expected discounted budget. Each producer faces a regular one-dimensional demand diffusion and chooses irreversible capacity expansion through a singular stochastic control problem. The regulator values weighted consumer surplus, producer profits, and an additional benefit from installed capacity, while accounting for the excess cost of raising public funds. We characterize producers’ optimal investment responses through free boundaries and reformulate the regulator’s stochastic optimization problem as a deterministic optimization over investment boundaries. An integral representation then recovers the optimal subsidy schedules that implement the resulting capacity policies. Under mild conditions on the diffusions, we establish existence of optimal schedules and uniqueness of the induced investment policies. Finally, we characterize the subsidy schedule and provide numerical illustrations.
\end{abstract}

\section{Introduction}
\label{sec:introduction}
\shortcites{BigernaEtAl2019,AtheyEtAl2022,CastilloEtAl2021,ZhangEtAl2024}

Investment subsidies influence both the timing and scale of productive
investment, with consequences for output, prices, and welfare.
\citet{McDonaldSiegel1986}, \citet{Pindyck1988}, and
\citet{DixitPindyck1994} explain how irreversibility and uncertainty
make the option to wait valuable. \citet{Guthrie2006} reviews how
regulatory incentives depend on the allocation of risk and the design
of policy. In a renewable-energy investment model,
\citet{BigernaEtAl2019} show that higher support can induce earlier but
smaller investment. Policymakers must weigh the welfare gains from encouraging investment
against the cost of financing subsidies. When public funds are limited
and markets are heterogeneous, this requires deciding both how to
design subsidies and how to allocate the available budget across markets.

We consider $N$ markets, each
served by a risk-neutral monopolist producing a non-durable good.
Demand in each market follows a regular one-dimensional diffusion.
Producers choose output and expand capacity irreversibly. Goods and
productive capacity are specific to each market, although demand shocks
may be correlated. The regulator commits to a
subsidy per unit of investment, maximizing weighted consumer
surplus, producer profits, and additional social capacity benefits.
The policy is subject to a common expected discounted subsidy budget
and also accounts for an excess cost of raising public funds.

Our main contribution is a method for constructing optimal subsidies from producers' investment boundaries, as established in Proposition~\ref{prop:regulatorCandidateVerification}. The key difficulty is that the restriction used in this construction does not cover every admissible nonnegative subsidy. To this end, Lemma~\ref{lem:regulatorPointwiseComparison} provides the comparison needed to establish optimality over the full admissible class. Under the stated assumptions, Proposition~\ref{prop:regulatorCandidateVerification} establishes existence of optimal subsidies and uniqueness of the induced investment policy from the given initial states. Theorem~\ref{thm:regulatorExistence} then determines the common budget multiplier and the optimal allocation across markets for every positive budget.

The additional welfare flow $\varepsilon_i A_i(C)$ represents benefits
of installed capacity that are not already included in the modeled
consumer surplus and producer profits. These may include preparedness
and infrastructure advantages from maintaining productive capabilities.
\citet{BaxterKing1993} allow public capital to enter utility, providing
an analogy for this separate welfare term.  Vaccine-capacity studies of \citet{AtheyEtAl2022} and
\citet{CastilloEtAl2021} illustrate how availability can generate
social gains exceeding manufacturers' commercial returns. We represent
such additional benefits through $A_i$, and use $\varepsilon_i$ to measure their importance to the regulator.

In a closely related monopoly-policy study,
\citet*{WenHagspielKort2017} study price support and investment-cost
reimbursement for a single irreversible capacity investment under
uncertainty. Their welfare-alignment results depend on the demand
specification and on whether support is conditional on an investment
trigger. Our setting allows continuing capacity expansion and jointly
chooses demand-dependent subsidy schedules across heterogeneous markets
sharing a budget. \citet{Pennings2000} studies subsidy packages financed
by profit taxation, while \citet{LukasThiergart2019} analyze the
interaction of debt financing, cash grants, and government bargaining.
\citet{ZhangEtAl2024} examine subsidy design under monopoly and
duopolistic competition, incorporating strategic investment interactions
within a market. Related investment models are developed
by \citet{BaldurssonKaratzas1997} and
\citet*{KardarasPavlisZervos2025} for competitive markets,
and by \citet{AlMotairiZervos2017} for a monopolistic structure.
In our setting, separate monopolistic markets are linked through
a common subsidy budget.

The allocation problem also connects to investment under limited
resources. \citet{Bank2005} studies monotone control under a dynamic
fuel constraint, while \citet*{ChiarollaFerrariRiedel2013} characterize
a social planner's allocation of irreversible investment among firms
under a time-dependent stochastic resource constraint. In our setting,
firms retain control of investment and the regulator acts through
subsidies. The aggregate expected expenditure constraint produces a
scalar budget multiplier, which supplements the excess financing cost
in the regulator's market-specific optimization problems.

The producer's expansion problem is a singular
stochastic control problem. \citet{KaratzasShreve1984} establish the equivalence
between a monotone-follower problem and an associated optimal
stopping problem, providing the foundation for analyzing investment
through marginal stopping decisions. \citet{LambertonZervos2013}
develop potential representations and distributional
variational-inequality characterizations for optimal stopping of
one-dimensional diffusions.

A related approach constructs investment policies through stochastic
representation. \citet{BankElKaroui2004} establish a representation
theorem for optional processes with applications to singular control
and allocation problems. Building on this approach,
\citet{RiedelSu2011} establish existence of optimal irreversible
investment under general uncertainty and characterize the policy
through a base-capacity process solving a backward
stochastic equation. For regular one-dimensional diffusions,
\citet{Ferrari2015} links the base-capacity process to a Markovian
investment boundary and derives an integral equation for that
boundary.

A closely related study is by
\citet*{DeAngelisFedericoFerrari2017}, who model market uncertainty and
unit investment costs as independent diffusions and allow a general
convex running cost. Our producer analysis characterizes irreversible
investment when marginal installation costs $\ell(x,c)$ depend jointly
on installed capacity and the same demand state that drives operating
profits. This specification allows demand to affect both investment
returns and the cost of expansion. The regulator's problem then
specializes the instrument to demand-dependent subsidies,
$\ell_i(x)=K_i-k_i(x)$, and uses the producer's boundary characterization
which allows us to explicitly construct the subsidy policy.

The resulting policy characterizations identify how fiscal scarcity
and additional capacity benefits shape support. Without additional capacity benefits, the constructed optimal schedules are nonincreasing in demand. When the budget is sufficiently large and
the excess cost of raising funds is zero, constant subsidies are optimal. Positive
capacity benefits can instead generate subsidies that initially
increase before declining. When the penalty on subsidy spending is high enough for a given market, subsidies stop once demand exceeds a certain level. The numerical
illustrations also show that larger per-unit subsidies need not imply
larger budget allocations: expenditure depends on the timing and scale
of the induced investment.

Section~\ref{sec:setup} presents the model,
Section~\ref{sec:producer} characterizes producers' investment responses,
and Section~\ref{sec:regulator} solves the regulator's problem.
Section~\ref{sec:illustrations} provides numerical illustrations.
The appendices contain the proofs.

\section{Setup: Separate Markets, Producers and Regulator}\label{sec:setup}

Consider the production of $N$ non-durable goods indexed by $i\in\{1,\ldots,N\}$, each good served by a
single risk-neutral producer. Markets are separate, and goods are non-substitutable: there is no trade, substitution in demand, or sharing of productive capacity across markets.
Producer $i$ chooses output $Q_i$ given capacity $C_i$ to serve its own market.

The profit function of the producer of the ith market per unit time is given by
\begin{equation}\label{eq:profit function}
 \Pi_i(P_i,Q_i)
 :=
 P_i(Q_i)Q_i
 -
 \frac{\beta_i\lambda_i^{-1/\beta_i}}
      {(1+\beta_i)^{1+1/\beta_i}}
 \frac{Q_i^{1+1/\beta_i}}{C_i^{\alpha_i/\beta_i}},
 \qquad
 \beta_i>0,\quad
 0<\alpha_i<1,\quad
 \lambda_i>0.
\end{equation} 

In addition, since each good is non-durable the market for each good satisfies the following market clearing condition
\begin{equation}\label{eq:market clearing}
Q_i(t)=\frac{D_{i}(t)}{P^{\delta_i}_{i}(t)}, \quad \delta_i>1,\quad i=1,...,N,    
\end{equation}
where $(D_i)_{i=1}^N$ represents the stochastic base preferences of an exogenous representative consumer. The base demand processes may be correlated and are unaffected by changes in production or
investment.
 
Each producer chooses its optimal output per unit time.
Maximizing \eqref{eq:profit function} with respect to $Q_i$,
taking into account the dependence of $P_i(Q_i)$ on output,
and using \eqref{eq:market clearing} gives:
\begin{equation}
 Q_i^*=\lambda_i(1+\beta_i)
      \left(1-\frac{1}{\delta_i}\right)^{\beta_i} C_i^{\alpha_i} P_i^{\beta_i},
 \qquad
 \Pi^*_i=\frac{\delta_i+\beta_i}{\delta_i(1+\beta_i)}P_iQ_i^*.
\end{equation}
For later use, put
\begin{equation}\label{eq:constants}
 \eta_i=\frac{\alpha_i(\delta_i-1)}{\delta_i+\beta_i},\qquad
 \chi_i=\frac{\delta_i(1+\beta_i)}{\delta_i+\beta_i}.
\end{equation}
Note that $\delta_i>1$ together with $0<\alpha_i<1$ and $\beta_i>0$ implies $0<\eta_i<1$.
For convenience we use the following normalization \begin{equation}\label{eq:normalizedDemand}
 X_{i}(t):=
 \left[\lambda_i(1+\beta_i)
 \left(1-\frac{1}{\delta_i}\right)^{\beta_i}
 \right]^{\frac{\delta_i-1}{\delta_i+\beta_i}}
 D_{i}^{\frac{1+\beta_i}{\delta_i+\beta_i}}(t).
\end{equation}
Consequently, the operating profit can be re-written as
\begin{equation}\label{eq:profit}
 \Pi_i(X_i,C_i)=\frac{X_iC_i^{\eta_i}}{\chi_i},
\end{equation}
and the market price is
\begin{equation}\label{eq:price}
 \mathsf P_i=
 \left[
 \frac{X_i}{\lambda_i (1+\beta_i)((\delta_i-1)/\delta_i)^{\beta_i} C_i^{\alpha_i-\eta_i}}
 \right]^\frac{1}{1+\beta_i}.
\end{equation}

On a complete probability space $(\Omega,\mathcal F,\mathbb P)$,
let $W=(W_1,\ldots,W_N)$ be an $N$-dimensional Brownian motion
with covariance matrix $\Gamma$ and $\Gamma_{ii}=1$.
For each market $i$, let
$\mathbb F_i=(\mathcal F_{it})_{t\ge0}$ be the completed,
right-continuous filtration generated by $W_i$.
Producer $i$'s decisions are adapted to $\mathbb F_i$.
We model  demand in market $i$ by
\[
 \mathrm dX_{i}(t)
 =\mu_i(X_{i}(t))\dd t+\sigma_i(X_{i}(t))\dd W_{i}(t),
 \qquad X_i(0)=x_i>0.
\]
Throughout the document we assume the following condition
\begin{assumption}\label{ass:diffusionDynamics}
The coefficients $\mu_i$ and $\sigma_i$ are locally Lipschitz
on $(0,\infty)$, with $\sigma_i>0$, and each $X_i$ is a regular,
nonexplosive diffusion that remains in $(0,\infty)$.
\end{assumption}

\subsection{Regulator and expansion policy}
Each producer can irreversibly expand its productive capacity at an
installation cost $K_i>0$ per unit. Producer $i$'s policy is adapted to the usual
augmented filtration generated by $W_i$. Initial capacity is
$C_i(0-)=c_i>0$, and future cash flows are discounted at the common rate
$r>0$.

To support this expansion, a regulator
commits to nonnegative demand-dependent subsidy
schedules $k=(k_1,\ldots,k_N)$. Producer $i$ takes $k_i$ as given and faces
the net investment cost per unit

\[
 \ell_i(x):=K_i-k_i(x).
\]

Taking its schedule as
given, each producer chooses investment and output to maximize
its own expected discounted profit. The regulator therefore evaluates
each schedule through the decisions it induces and their
effects on installed capacity, production, prices, and consumer
surplus.

Consumer surplus is total willingness to pay minus expenditure.
It equals
\begin{equation}
 \int_0^{Q_i}P_i(q)\dd q-P_i(Q_i)Q_i
 =\frac{P_i(Q_i)Q_i}{\delta_i-1},
\end{equation}

where we have used \eqref{eq:market clearing} to perform the integration.

At the producer's optimal output,
$P_i(Q_i^*)Q_i^*=xc^{\eta_i}$.
Thus instantaneous consumer surplus and operating profit are
\begin{equation}\label{eq:consumerSurplusFlow}
 \operatorname{\mathcal{CS}}_i(x,c)=\frac{xc^{\eta_i}}{\delta_i-1},
 \qquad
 \Pi_i(x,c)=\frac{xc^{\eta_i}}{\chi_i}.
\end{equation}
Given the regulator's schedules $k=(k_1,\ldots,k_N)$, suppose
for now that each producer has an optimal capacity response,
denoted by $C_i^{k_i}(t)$. Define the corresponding price by
\[
 P_i^{k_i}(t)=\mathsf P_i(X_i(t),C_i^{k_i}(t)).
\]
Using the common rate $r$, define expected
consumer surplus, company profit, and subsidy expenditure by
\begin{align}
 \mathcal{CS}_i(k_i)
 &=\E_{x_i,c_i}\!\left[
   \int_0^\infty e^{-rt}
   \frac{X_i(t)(C_i^{k_i}(t))^{\eta_i}}{\delta_i-1}
   \dd t\right],
   \label{eq:discountedConsumerSurplus}\\
 V_i^{k_i}(x_i,c_i)
 &=\E_{x_i,c_i}\!\bigg[
   \int_0^\infty e^{-rt}
   \frac{X_i(t)(C_i^{k_i}(t))^{\eta_i}}{\chi_i}
   \dd t-\int_{[0,\infty)}e^{-rt}
   \bigl[K_i-k_i(X_i(t))\bigr]
   \dd C_i^{k_i}(t)\bigg],
   \label{eq:regulatorProducerProfit}\\
 S_i(k_i)
 &=\E_{x_i,c_i}\!\left[
   \int_{[0,\infty)}e^{-rt}
   k_i(X_i(t))\dd C_i^{k_i}(t)
   \right].
   \label{eq:regulatorExpenditure}
\end{align}
The regulator also assigns a direct policy value to installed capacity.
For each market, let $A_i:[c_i,\infty)\to[0,\infty)$, normalized by
$A_i(c_i)=0$, and let $\varepsilon_i\ge0$ be its weight. Define
\begin{equation}\label{reg:capacityBenefitPayoff}
 \mathcal B_i(k_i):=\varepsilon_i\E_{x_i,c_i}\!\left[
   \int_0^\infty e^{-rt}A_i(C_i^{k_i}(t))\dd t\right].
\end{equation}
This reflects the public benefit of advanced infrastructure.

Let $\vartheta>0$ measure the relative weight on consumer surplus.
The regulator's objective combines market welfare and capacity benefits,
net of public expenditure:
\begin{equation}\label{eq:regulatorPayoff}
 \begin{aligned}
 J_{\vartheta}(k)=\frac{1}{1+\vartheta}\Bigg[
   &\vartheta\sum_{i=1}^N\mathcal{CS}_i(k_i)
   +\sum_{i=1}^N\bigl(V_i^{k_i}(x_i,c_i)+\mathcal B_i(k_i)\bigr)-(1+\zeta)S(k)\Bigg],\qquad \zeta\ge0,
 \end{aligned}
\end{equation}

where $S(k):=\sum_{i=1}^N S_i(k_i)$
Here $\zeta$ is the excess resource cost of raising one unit of public
funds. In particular, $\zeta=0$ treats subsidies as costless transfers.

Given a common budget $B\ge0$, the regulator chooses the schedules jointly:
\begin{equation}\label{eq:regulatorProblem}
 V_{\mathrm{reg}}(B)=
 \sup_{\substack{k\in\mathcal K_+\\S(k)\le B}}J_{\vartheta}(k).
\end{equation}
Here $\mathcal K_+$ denotes the admissible schedules satisfying
the producer's conditions imposed below.

We fix a producer and omit its index. The assumptions
apply separately to each market.

\section{The Producer's Control Problem}\label{sec:producer}

The producer's state space is
\[
 \mathcal D:=\{(x,c)\in\mathbb R^2:x>0,\ c>0\}.
\]

Before studying the producer's problem, we introduce the diffusion
quantities used below \cite{BorodinSalminen2002}.
The infinitesimal generator of $X$ is
\[
 \mathcal L f(x)=\tfrac12\sigma^2(x)f''(x)+\mu(x)f'(x).
\]
For $\tau_y:=\inf\{t\ge0:X(t)\ge y\}$,
\begin{equation}\label{eq:harmonic}
 (\mathcal L-r)\psi=0,\qquad
 \E_x[e^{-r\tau_y};\tau_y<\infty]
 =\frac{\psi(x)}{\psi(y)},\qquad 0<x\le y,
\end{equation}
where $\psi$ denotes the positive increasing solution of
$(\mathcal L-r)\psi=0$, characterized by the hitting-time identity
in \eqref{eq:harmonic}.

The expected discounted demand is
\begin{equation}\label{eq:resolvent}
 U(x):=\E_x\!\left[\int_0^\infty e^{-rt}X(t)\dd t\right],
 \qquad (\mathcal L-r)U(x)=-x.
\end{equation}
We also define
\begin{equation}\label{eq:qH}
 q(x):=\frac{\psi(x)}{\psi'(x)},\qquad
 H(x):=U(x)-q(x)U'(x).
\end{equation}

\begin{assumption}\label{ass:diffusion}
The functions $U$ and $\psi$ belong to
$C^2((0,\infty))$, with $U<\infty$, $\psi>0$, and $\psi'>0$.
Moreover $ \lim_{x\to\infty}U(x)=\infty$.
\end{assumption}

\begin{corollary}\label{cor:resolventProperties}
If
$U(x)<\infty$ for every $x>0$, then
\[
\lim_{x\to\infty}\frac{U(x)}{\psi(x)}=0,
\qquad
0<H(x)<\frac{x}{r},
\qquad
H'(x)
=\frac{2q(x)}{\sigma^2(x)}
 \bigl[x-rH(x)\bigr]>0.
\]
In particular, $H(0+)=0$, $\lim_{x\rightarrow\infty}H(x)=\infty$ and $U/\psi$ is strictly
decreasing.
\end{corollary}

The producer chooses an adapted, right-continuous, nondecreasing
capacity process with $C_{0-}=c>0$, finite on bounded time intervals.
Expenditure includes the initial jump and is integrated along each
installation path:
\begin{equation}\label{eq:jumpConvention}
 \int_{[0,T]}e^{-rt}\ell(X(t),C(t))\dd C_t
 :={}\int_0^T e^{-rt}\ell(X(t),C(t))\dd C^{\rm c}(t)\notag+\sum_{0\le t\le T}e^{-rt}
       \int_{C(t-)}^{C(t)}\ell(X(t),u)\dd u,
\end{equation}
where $C^c$ denotes the continuous part of the process.

\begin{remark}
In this section we allow for general marginal installation cost $\ell(x,c)$.    
\end{remark}

Note that during continuous capacity changes,
\begin{equation}\label{eq:priceDynamics}
 \frac{\mathrm dP(t)}{P(t)}
 =\left[\frac{\mu(X(t))}{(1+\beta)X(t)}
       -\frac{\beta\sigma^2(X(t))}{2(1+\beta)^2X(t)^2}\right]\dd t
  +\frac{\sigma(X(t))}{(1+\beta)X(t)}\dd W_t
  -\frac{\alpha}{\delta+\beta}\frac{\mathrm dC(t)}{C(t)}.
\end{equation}

\begin{definition}
The set $\mathcal A$ of admissible expansion strategies consists of
all adapted, c\`adl\`ag, nondecreasing capacity processes $C$ with
$C_{0-}=c$ that satisfy
\begin{equation}\label{eq:admissibility}
 \E_{x,c}\!\left[
 \int_{[0,\infty)}e^{-rt}\ell(X(t),C(t))\dd C(t)
 \right]<\infty,
\end{equation}
\end{definition}
The producer's value function is
\begin{equation}\label{eq:valueProblem}
 V(x,c)=\sup_{C\in\mathcal A}
 \E_{x,c}\!\left[
 \int_0^\infty e^{-rt}\frac{X(t)C^\eta(t)}{\chi}\dd t
 -\int_{[0,\infty)}e^{-rt}\ell(X(t),C(t))\dd C(t)
 \right].
\end{equation}
The associated Hamilton--Jacobi--Bellman equation is
\begin{equation}\label{eq:HJB}
 \max\left\{
 (\mathcal L-r)w(x,c)+\frac{xc^\eta}{\chi},\quad
 w_c(x,c)-\ell(x,c)
 \right\}=0,  \quad(x,c)\in\mathcal D.
\end{equation}

We proceed in three steps: construct the boundary and candidate value,
prove that the candidate satisfies \eqref{eq:HJB} in a measure sense,
and verify optimality of the control.

\subsection{Boundary and candidate construction}
For the fixed subsidy schedule, we seek a continuous, strictly
increasing demand boundary $b:(0,\infty)\rightarrow(0,\infty)$. 

The demand boundary separates the state space into the waiting
and investment regions
\begin{equation}\label{eq:producerRegions}
 \mathcal W:=\{(x,c)\in\mathcal D:x<b(c)\},
 \qquad
 \mathcal I:=\{(x,c)\in\mathcal D:x\ge b(c)\}.
\end{equation}
The proposed policy keeps capacity fixed in $\mathcal W$.
If $x>b(c)$, it increases capacity immediately until the
boundary is reached. Existence of the boundary and optimality of this policy are
established below.

For a fixed capacity level, define the net value of installing one
marginal unit immediately by
\begin{equation}\label{eq:marginalGain}
 g(x,c)=\frac{\eta}{\chi}c^{\eta-1}U(x)-\ell(x,c).
\end{equation} 
In the
waiting region, the candidate solution has the form
\begin{equation}\label{eq:waitingValue}
 w(x,c)=\frac{c^\eta}{\chi}U(x)+F(c)\psi(x).
\end{equation}

To determine $F$ and the boundary $b$, we appeal to the so-called  "smooth-fit principle"
condition of singular control. In particular, using \eqref{eq:waitingValue}, this gives
\begin{align}
 \lim_{x\uparrow b(c)}w_c(x,c)
 &=\frac{\eta}{\chi}c^{\eta-1}U(b(c))
   +F'(c)\psi(b(c))=\ell(b(c),c),
 \label{eq:marginalFit}\\
 \lim_{x\uparrow b(c)}w_{cx}(x,c)
 &=\frac{\eta}{\chi}c^{\eta-1}U'(b(c))
   +F'(c)\psi'(b(c))=\ell_x(b(c),c).
 \label{eq:smoothFit}
\end{align}
Multiplying \eqref{eq:smoothFit} by $q(b(c))$ and subtracting
it from \eqref{eq:marginalFit} eliminates $F'(c)$ and yields
\begin{equation}\label{eq:boundaryMain}
 \ell(b(c),c)
 -q(b(c))\ell_x(b(c),c)
 =\frac{\eta}{\chi}c^{\eta-1}H(b(c)).
\end{equation}

We now give conditions that make this equation define the boundary and
verify optimality. All the functions used in these conditions are determined
by the diffusion and $k$. Put
\begin{align}
 \Theta(x,c)
   &:=\ell(x,c)
 -q(x)\ell_x(x,c),
      \label{eq:adjustedCost}\\
 \mathcal T(x,c)
   &:=\frac{c^{1-\eta}\Theta(x,c)}{H(x)}.
       \label{eq:thresholdFunction}
\end{align}
The boundary equation is simply
\begin{equation}\label{eq:boundaryT}
 \mathcal T(b(c),c)=\frac{\eta}{\chi}.
\end{equation}

\begin{assumption}\label{ass:subsidy}
The net installation cost $0<\underline{\ell}<\ell\leq K$ belongs to $C^1(\mathcal D)$.
The continuous function $\mathcal T$ satisfies
\begin{equation}\label{eq:thresholdSigns}
 \begin{gathered}
 x\longmapsto\mathcal T(x,c)\text{ is strictly decreasing for each }c>0,\\
 c\longmapsto\mathcal T(x,c)\text{ is strictly increasing for each }x>0.
 \end{gathered}
\end{equation}
Its endpoint limits are
\begin{align}
 \lim_{x\downarrow0}\mathcal T(x,c)&=\infty,&
 \lim_{x\to\infty}\mathcal T(x,c)&=0
 &&\text{for each }c>0,\label{eq:thresholdXLimits}\\
 \lim_{c\downarrow0}\mathcal T(x,c)&=0,&
 \lim_{c\to\infty}\mathcal T(x,c)&=\infty
 &&\text{for each }x>0.\label{eq:thresholdCLimits}
\end{align}
\end{assumption}

\begin{proposition}\label{prop:boundary}
Under Assumptions~\ref{ass:diffusion} and \ref{ass:subsidy}, equation
\eqref{eq:boundaryT} has a unique solution $b(c)>0$ for every $c>0$.
The map $b:(0,\infty)\to(0,\infty)$ is a strictly increasing continuous
bijection. For each $c$, the function $x\mapsto g(x,c)/\psi(x)$
has its unique global maximum at $b(c)$. Also, $g(x,c)>0$ for $x\ge b(c)$.

\end{proposition}

Let $\Phi:(0,\infty)\rightarrow (0,\infty)$ be defined by $\Phi=b^{-1}$.  By
Proposition~\ref{prop:boundary}, $\Phi$ is continuous and strictly increasing. Solving \eqref{eq:boundaryT} in terms of $\Phi$ gives
\begin{equation}\label{eq:producerCapacitySchedule}
\Theta(y,\Phi(y)) \Phi^{1-\eta}(y)=
 \frac{\eta H(y)}{\chi}
\end{equation}

\begin{remark}
On any open set where $\mathcal T$ is continuously differentiable and 
$\partial_c\mathcal T(b(c),c)>0>\partial_x\mathcal T(b(c),c)$,  $b(c)$
is $C^1$ and
\begin{equation}\label{eq:boundaryDerivative}
 b'(c)=-\frac{\partial_c\mathcal T(b(c),c)}
                  {\partial_x\mathcal T(b(c),c)}>0.
\end{equation}
\end{remark}

\begin{assumption}\label{ass:tail}
For every $c>0$, the boundary in Proposition~\ref{prop:boundary} satisfies
\begin{equation}\label{eq:tail}
 \int_c^\infty
 u^{\eta-1}\frac{U(b(u))}{\psi(b(u))}\dd u<\infty.
\end{equation}
\end{assumption}

Equation \eqref{eq:marginalFit} determines $F$ up to an additive constant:
\[
 F'(c)=-\frac{g(b(c),c)}{\psi(b(c))}.
\]
We select the coefficient  via the asymptotic condition
$\lim_{c\to\infty}F(c)=0$. Thus,
\begin{equation}\label{eq:F}
 F(c)=\int_c^\infty
 \frac{\frac{\eta}{\chi}u^{\eta-1}U(b(u))
                -\ell(b(u),u)}{\psi(b(u))}\dd u.
\end{equation}
Since
$0<g(b(u),u)<(\eta/\chi)u^{\eta-1}U(b(u))$,
Assumption~\ref{ass:tail} ensures that $F$ is finite and nonnegative.

Thus, the candidate value on the whole state space is given by
\begin{equation}\label{eq:candidate}
  w(x,c)=
 \begin{cases}
 \displaystyle\frac{c^\eta}{\chi}U(x)+F(c)\psi(x),
                         &x\le b(c),\\[6pt]
 \displaystyle\frac{\Phi(x)^\eta}{\chi}U(x)+F(\Phi(x))\psi(x)
                -\int_c^{c\vee\Phi(x)}\ell(x,u)\dd u,
                         &x>b(c)
 \end{cases}
\end{equation}
The second branch includes the initial adjustment from $c$ to $c\vee\Phi(x)$ and
pays the installation cost along that capacity path. Defining
\begin{equation}\label{eq:ConeMarginalValue}
 \mathcal{Z}(x,c):=
 \begin{cases}
 \displaystyle
 \frac{\psi(x)}{\psi(b(c))}g(b(c),c),
                                      &x<b(c),\\[4pt]
 g(x,c),                              &x\ge b(c).
 \end{cases}
\end{equation}

and after integration by part we get the equivalent
expression
\begin{equation}\label{eq:candidateLayers}
  w(x,c)=\frac{c^\eta}{\chi}U(x)
       +\int_c^{\infty}\mathcal{Z}(x,u)\dd u
\end{equation}


We now move to establish that the candidate $w$ satisfies the HJB equation.

Since $k$ is only assumed to be $C^1$, we interpret the diffusion term as a measure. For $h\in C^1((0,\infty))$ with $h'$ of
locally bounded variation, define
\begin{equation}\label{eq:measureGenerator}
 \mathcal G_r h(\mathrm dx)
 :=\frac12\sigma^2(x)\,\mathrm d h'(x)
   +[\mu(x)h'(x)-rh(x)]\,\mathrm dx.
\end{equation}
Here $\mathrm d h'$ is the distributional second derivative of $h$.
If $h$ is $C^2$, then
$\mathcal G_r h=[(\mathcal L-r)h](x)\,\mathrm dx$.

\begin{lemma}\label{lem:HJB}
Suppose that Assumptions~\ref{ass:diffusion}, \ref{ass:subsidy},
and \ref{ass:tail} hold. The candidate $w$ in
\eqref{eq:candidate} is 
$C^1$, with $w_x(\cdot,c)$ locally of bounded
variation for every $c>0$. It satisfies the HJB equation \eqref{eq:HJB} 
\end{lemma}

\subsection{Verification of the optimal investment policy}

The next result proves optimality of the construction. The proof uses
one-dimensional marginal stopping problems, so it does not require
a classical two-variable It\^o formula for $w$.

\begin{theorem}\label{thm:verification}
Suppose that Assumptions~\ref{ass:diffusion}--\ref{ass:tail}
hold. Then the value function $V$ of
\eqref{eq:valueProblem} equals the candidate $w$ in \eqref{eq:candidate}. The unique optimal control process is
\begin{equation}\label{eq:reflection}
 C_t^*=c\vee\Phi(\overline X(t)),\qquad \overline X(t)=\sup_{0\le s\le t}X_s,
 \qquad C^*(0-)=c.
\end{equation}
\end{theorem}

\section{The Regulator's Problem}
\label{sec:regulator}

For the regulator's problem the policy instrument is a subsidy indexed by demand so
$\ell(x)=K-k(x)$ and $\Theta(x)=\ell(x)-q(x)\ell_x(x)$. 

\subsection{Assumptions and admissible schedules}

Suppress the market index temporarily and retain the diffusion quantities
$U$, $\psi$, $q$, and $H$ in
\eqref{eq:harmonic}--\eqref{eq:qH}. Define
\begin{equation}\label{reg:stateWeight}
 \omega(y):=\frac{\psi'(y)}{\psi^2(y)}.
\end{equation}

\begin{assumption}\label{ass:regulatorExistence}
Assumption~\ref{ass:diffusion} holds, and for every $z>0$,
\begin{equation}\label{reg:integrability}
 \int_z^\infty \omega(y)H^{1/(1-\eta)}(y)\,\dd y<\infty.
\end{equation}
\end{assumption}

We impose the following conditions on how the regulator socially values additional capacity

\begin{assumption}\label{ass:capacityBenefit}
Let
\[
A(z)=\int_c^z \frac{G_A(u)}{u}\,\dd u,
\qquad z\ge c,
\]
where \(G_A\in C^2([c,\infty))\) satisfies
\[
G_A(c)=0,\qquad G_A(z)>0\quad\text{for }z>c,~~
G_A'(z)\ge0,~~ G_A''(z)\le0
\quad\text{for }z\ge c,~~
\lim_{z\to\infty}\frac{G_A(z)}{z}=0.
\]
Assume also that
\begin{equation}\label{reg:capacityBenefitStrength}
\frac{\varepsilon}{r}G_A'(c)<(1-\eta)K.
\end{equation}
\end{assumption}
Thus $G_A=zA'$ is increasing and concave, $A'(c)=0$, and $A'(z)\to0$.
 Also, $0\le A'(z)\le G'_A(c)$.

Assumptions~\ref{ass:regulatorExistence} and
\ref{ass:capacityBenefit} are imposed for every market throughout the
regulator's construction.

We restore the market index $i$ and fix the initial states $(x_i,c_i)$. The admissibility class is
\begin{equation}\label{reg:stateAdmissibility}
\mathcal K_{i,+}
:=
\left\{
\begin{array}{l}
k\in C^1((0,\infty)):\ k\ge0, S_i(k)<\infty \\
\text{Assumptions~\ref{ass:subsidy} and~\ref{ass:tail} hold}
\end{array}
\right\},\qquad \mathcal K_+:=\prod_i\mathcal K_{i,+}.
\end{equation}

Under Assumption~\ref{ass:regulatorExistence}, Theorem~\ref{thm:verification}
applies to every $k\in\mathcal K_{i,+}$, since the definition of
$\mathcal K_{i,+}$ ensures Assumptions~\ref{ass:subsidy}
and~\ref{ass:tail}.

At this point we define the following quantities
\begin{equation}
\begin{aligned}\label{eq:YI}
Y_i(k_i)&:=\E_{x_i,c_i}\!\left[
\int_0^\infty e^{-rt}X_i(t)(C_i^{k_i}(t))^{\eta_i}\dd t\right],
& I_i(k_i)&:=\E_{x_i,c_i}\!\left[
\int_{[0,\infty)}e^{-rt}\dd C_i^{k_i}(t)\right]
\end{aligned}
\end{equation}

Cancelling subsidy receipts against public payments in
\eqref{eq:regulatorPayoff}, put
\begin{equation}\label{reg:welfareReduction}
 h_i:=\frac1{\chi_i}+\frac{\vartheta}{\delta_i-1}
\end{equation}
Thus, the regulator solves
\begin{equation}\label{reg:constrainedProblem}
 V_{\rm reg}(B)=\frac1{1+\vartheta}
 \sup_{\substack{k\in\mathcal K_+\\\sum_iS_i(k_i)\le B}}
 \sum_i[h_iY_i(k_i)-K_iI_i(k_i)-\zeta S_i(k_i)+\mathcal B_i(k_i)].
\end{equation}
For a multiplier $\Lambda\ge0$, let $\rho:=\zeta+\Lambda$ and define
\begin{equation}\label{reg:penalizedMarket}
 \mathcal J_{i,\rho}(k_i):=h_iY_i(k_i)-K_iI_i(k_i)-\rho S_i(k_i)+\mathcal B_i(k_i).
\end{equation}
Thus, the common expected budget couples the market-level policy choices.
\begin{proposition}\label{prop:regulatorVerification}
Suppose $k_i^*$ maximizes $\mathcal J_{i,\zeta+\Lambda^*}$ over
$\mathcal K_{i,+}$ for every market, where $\Lambda^*\ge0$, and
\begin{equation}\label{reg:budgetComplementarity}
 S(k^*)\le B,\qquad \Lambda^*[B-S(k^*)]=0.
\end{equation}
Then $k^*$ solves \eqref{reg:constrainedProblem}, and
\begin{equation}\label{reg:dualValue}
 V_{\rm reg}(B)=\frac{\Lambda^*B+
                  \sum_i\mathcal J_{i,\zeta+\Lambda^*}(k_i^*)}
                  {1+\vartheta}.
\end{equation}
\end{proposition}
\begin{proof}
For every feasible $k$, market optimality implies
\[
 (1+\vartheta)J_\vartheta(k)
 =\sum_i\mathcal J_{i,\zeta+\Lambda^*}(k_i)+\Lambda^*S(k)
 \le\sum_i\mathcal J_{i,\zeta+\Lambda^*}(k_i^*)+\Lambda^*B
 =(1+\vartheta)J_\vartheta(k^*).
\]
The last equality is complementary slackness.
\end{proof}

Applying the hitting-time identity to \eqref{eq:regulatorExpenditure}, \eqref{reg:capacityBenefitPayoff} and \eqref{eq:YI} gives
\begin{align}
 Y_i(k_i)
 &=(c_i\vee \Phi_i(x_i))^{\eta_i}U_{i}(x_i)
 +\psi_{i}(x_i)\int_{c_i\vee \Phi_i(x_i)}^\infty
   \frac{\eta_i u^{\eta_i-1}U_{i}(b_i(u))}
        {\psi_{i}(b_i(u))}\dd u,
 \label{reg:Yboundary}\\
 I_i(k_i)
 &=c_i\vee \Phi_i(x_i)-c_i+\psi_{i}(x_i)\int_{c_i\vee \Phi_i(x_i)}^\infty
   \frac{\dd u}{\psi_{i}(b_i(u))},
 \label{reg:Iboundary}\\
 S_i(k_i)
 &=k_i(x_i)(c_i\vee \Phi_i(x_i)-c_i)+\psi_{i}(x_i)\int_{c_i\vee \Phi_i(x_i)}^\infty
   \frac{k_i(b_i(u))}{\psi_{i}(b_i(u))}\dd u,\\
   \mathcal B_i(k_i)&=\frac{\varepsilon_i}{r}A_i\!\left(c_i\vee\Phi_i(x_i)\right)+\frac{\varepsilon_i}{r}\psi_{i}(x_i)\int_{c_i\vee\Phi_i(x_i)}^\infty\frac{A_i'(u)}{\psi_{i}(b_i(u))}\,du.
 \label{reg:Sboundary}
\end{align}

\begin{remark}
For fixed marginal demand dynamics and initial states, the optimal subsidies and regulator's value are independent of $\Gamma$. Each producer's problem uses only its own demand process and filtration but correlation can affect the distribution of aggregate expenditure.
\end{remark}

\subsection{The Regulator's Solution}
\label{sec:regulatorSolution}

The regulator chooses a subsidy schedule $k$, anticipating the
optimal response of each producer. We formulate her choice in terms of the capacity
boundaries $\lbrace\Phi_i\rbrace_{i=1}^N$ she wishes to induce, restricting attention to boundaries
implementable through admissible nonnegative subsidies satisfying the
producer's boundary condition \eqref{eq:boundaryT}. Since the regulator can solve each market separately then for a
market $i$ and fixed $\rho$, we rewrite the regulator's objective in
\eqref{reg:penalizedMarket} in terms of $\Phi_i$ and construct a candidate
boundary by pointwise maximization. We then recover its implementing
subsidy from the boundary condition and verify its admissibility,
ensuring that the desired investment response is optimal for the
producer. Finally, we determine the common budget multiplier that links
the market schedules through the aggregate budget constraint.

Fix a market and suppress its index. Throughout this part, $(x,c)$ are the given
initial states, $y$ is a variable demand state, and $z$ is a proposed capacity. Define
\begin{equation}\label{reg:adjustedMarginalReturn}
 M(y,z):=\frac{\eta}{\chi}H(y)z^{\eta-1},\qquad
 M_z(y,z)=-\frac{1-\eta}{z}M(y,z)<0.
\end{equation}
The free boundary equation is simply $M(y,\Phi(y))=\Theta(y)$.

For an admissible schedule $k$, every capacity level $u>c$ is
installed at time zero if $u\le\Phi(x)$ and otherwise at the first hitting time of $b(u)$.
The assumptions imply $\psi(y)\to\infty$. Since
$(U/\psi)'=-\omega H$ and $U/\psi\to0$,
\begin{equation}\label{reg:revenueTailIdentity}
 \frac{U(v)}{\psi(v)}=\int_v^\infty\omega(y)H(y)\dd y,
 \qquad
 \frac1{\psi(v)}=\int_v^\infty\omega(y)\dd y.
\end{equation}
Substituting these identities into \eqref{reg:Yboundary}
gives 
\[
\begin{aligned}
Y(k)-c^\eta U(x)
&=\eta\psi(x)\Biggl[
\int_c^{c\vee\Phi(x)}u^{\eta-1}
\int_x^\infty\omega(y)H(y)\dd y\dd u+\int_{c\vee\Phi(x)}^\infty u^{\eta-1}
\int_{b(u)}^\infty\omega(y)H(y)\dd y\dd u
\Biggr].
\end{aligned}
\]
Using the fact that $
\{c<u\le c\vee\Phi(x),\ y\ge x\}\cup\{u>c\vee\Phi(x),\ y\ge b(u)\}=\{y\ge x,\ c<u\le\Phi(y)\}$ and Tonelli's theorem we arrive at
\begin{equation}
 Y(k)=c^\eta U(x)+\psi(x)\int_x^\infty
       \omega(y)H(y)[\Phi^\eta(y)-c^\eta]_+\dd y,
       \label{reg:stateRevenue}
\end{equation}

Similar calculations give
\begin{equation}
 I(k)=\psi(x)\int_x^\infty
       \omega(y)[\Phi(y)-c]_+\dd y.
       \label{reg:stateInvestment}
\end{equation}

Finally, the same identities give
\begin{equation}\label{reg:stateCapacityBenefit}
 \mathcal B(k)=\frac{\varepsilon}{r}\psi(x)\int_x^\infty
       \omega(y)A(c\vee\Phi(y))\dd y.
\end{equation}
Since $A(z)\le G'_A(c)(z-c)$ for $z\ge c$, this payoff is finite whenever
discounted investment is finite.

For the net cost, \eqref{eq:adjustedCost} gives
\begin{equation}\label{reg:netCostTailIdentity}
 \left(\frac{K-k(y)}{\psi(y)}\right)'
       =-\omega(y)\Theta(y),\qquad
 \frac{K-k(v)}{\psi(v)}
       =\int_v^\infty\omega(y)\Theta(y)\dd y.
\end{equation}
The terminal term vanishes because $k$ is bounded. Since $\Theta(y)>0$,
Tonelli again applies and gives
\begin{equation}\label{reg:stateNetCost}
 KI(k)-S(k)=\psi(x)\int_x^\infty
       \omega(y)\Theta(y)[\Phi(y)-c]_+\dd y.
\end{equation}
Finally, the subsidy cost is written as
\begin{equation}\label{reg:boundaryTransferRepresentation}
 S(k)=\psi(x)\int_x^\infty
       \omega(y)[K-M(y,\Phi(y))][\Phi(y)-c]_+\dd y.
\end{equation}

Substituting these representations into \eqref{reg:penalizedMarket} gives
\begin{equation}\label{reg:boundaryObjectiveRepresentation}
 \mathcal J_\rho(k)=hc^\eta U(x)+\psi(x)\int_x^\infty
       \omega(y)F_\rho(y,\Phi(y))\dd y,
\end{equation}
where, for $z>c$,
\begin{equation}\label{reg:boundaryScalarObjective}
 F_\rho(y,z):=hH(y)(z^\eta-c^\eta)
       -(1+\rho)K(z-c)+\rho M(y,z)(z-c)+\frac{\varepsilon}{r}A(z),
\end{equation}
and $F_\rho(y,z)=0$ for $0<z\le c$. 
Finally, subtracting \eqref{reg:netCostTailIdentity} from $K$ times
the second identity in \eqref{reg:revenueTailIdentity} shows that
every bounded implementing schedule satisfies
\begin{equation}\label{reg:boundaryScheduleRecovery}
 k(y)=\psi(y)\int_y^\infty
       [K-M(t,\Phi(t))]\omega(t)\dd t.
\end{equation}
Note that, we also need to ensure that the capacity boundary generated by the regulator's pointwise maximization will result in a nonnegative subsidy schedule $k(y)\geq 0$.

For $z>c$, direct differentiation gives
\begin{equation}\label{reg:scalarObjectiveDerivative}
 \partial_zF_\rho(y,z)
 =M(y,z)\left[h\chi+\rho\eta+
                   \rho(1-\eta)\frac cz\right]+\frac{\varepsilon}{r}A'(z)-(1+\rho)K.
\end{equation}

Also define
\begin{equation}\label{reg:entryAdjustedReturn}
 d_\rho:=\frac{(1+\rho)K}{h\chi+\rho},\qquad
 y^{\text{entry}}_\rho:=H^{-1}\!\left(\frac{\chi d_\rho}{\eta}c^{1-\eta}\right),
\end{equation}

where $d_\rho$ is the critical value of $M(y,c)$ at which capacity expansion from $c$ becomes worthwhile, and $y^{\text{entry}}_\rho$ is the corresponding entry demand threshold, determined by $M(y^{\text{entry}}_\rho,c)=d_\rho$.

To determine the optimal capacity above $c$, set
$\partial_zF_\rho(y,z)=0$ and substitute the definition of $M$ from
\eqref{reg:adjustedMarginalReturn}. Rearranging gives
\begin{equation}\label{reg:scalarThresholdMap}
 z^{1-\eta}\left[(1+\rho)K-\frac{\varepsilon}{r}A'(z)\right]
 =\frac{\eta}{\chi}H(y)
 \left[h\chi+\rho\eta+\rho(1-\eta)\frac cz\right],\qquad z\ge c.
\end{equation}
For fixed $y$, the left-hand side is strictly increasing to infinity,
while the right-hand side is positive and nonincreasing in $z$. These properties are shown in
Appendix~\ref{app:regulatorExistenceProof}.
At $z=c$, equality holds exactly when $y=y^{\text{entry}}_\rho$, since $A'(c)=0$.
If $y\le y^{\text{entry}}_\rho$, then $\partial_zF_\rho(y,z)<0$ for every $z>c$,
so every $z\le c$ maximizes $F_\rho(y,\cdot)$.
If $y>y^{\text{entry}}_\rho$, there is a unique maximizer $\Phi_\rho^{\rm rel}(y)>c$ which is characterized directly by the first-order condition
\begin{equation}\label{reg:regulatorFOC}
 M(y,\Phi_\rho^{\rm rel}(y))
 \left[h\chi+\rho\eta+\rho(1-\eta)\frac{c}{\Phi_\rho^{\rm rel}(y)}\right]
 +\frac{\varepsilon}{r}A'(\Phi_\rho^{\rm rel}(y))=(1+\rho)K.
\end{equation}

Note that this solution is not necessarily admissible at this point.

 We define the no-subsidy boundary capacity from $M(y,\Phi_K(y))=K$ and  \eqref{reg:adjustedMarginalReturn}
\begin{equation}\label{reg:noSubsidyCapacity}
 \Phi_K(y):=\left(\frac{\eta H(y)}{\chi K}\right)^{1/(1-\eta)}.
\end{equation}
The producer’s capacity boundary induced by the regulator is
\begin{equation}\label{reg:candidateBoundary}
 \Phi_\rho(y):=
 \begin{cases}
  \left(\dfrac{\eta H(y)}{\chi d_\rho}\right)^{1/(1-\eta)},
                                      &0<y\le y^{\text{entry}}_\rho,\\[5pt]
  \Phi_\rho^{\rm rel}(y)\vee\Phi_K(y),&y>y^{\text{entry}}_\rho,
 \end{cases}
 \qquad m_\rho(y):=M(y,\Phi_\rho(y)).
\end{equation}
The first branch satisfies $m_\rho(y)=d_\rho$ and completes the boundary below the expansion level while the
second enforces $m_\rho\le K$ by implementing a larger producer's boundary, $\Phi_K(y)$ whenever $M(y,\Phi^{\rm rel}_\rho(y))>K$. The implementing subsidy is
\begin{equation}\label{reg:candidateOptimalSchedule}
 k_\rho(y)=\psi(y)\int_y^\infty[K-m_\rho(t)]\omega(t)\dd t,
\end{equation}
which by definition is nonnegative.

\begin{figure}[htbp]
\centering
\includegraphics[width=0.63\textwidth]{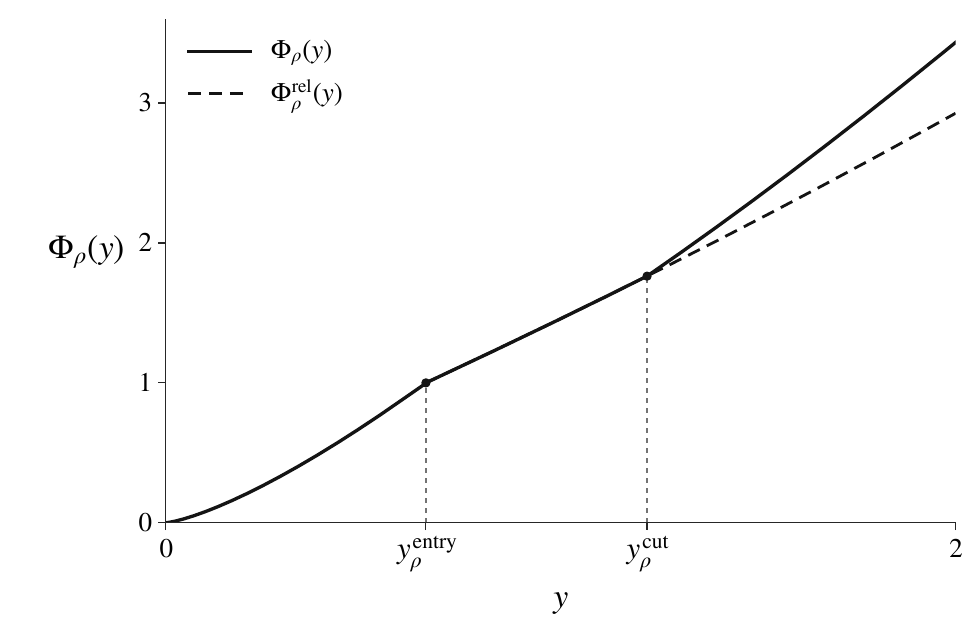}
\caption{Illustration of the optimal capacity boundary $\Phi_\rho(y)$ (solid) and
relaxed boundary $\Phi_\rho^{\rm rel}(y)$ (dashed).
Vertical lines mark entry and the subsidy cutoff.}
\label{fig:optimalCapacity}
\end{figure}

\begin{proposition}\label{prop:regulatorCandidateVerification}
For every $\rho\ge0$ the candidate
$\Phi_\rho$ is a continuous strictly increasing bijection.
The schedule $k_\rho$ in \eqref{reg:candidateOptimalSchedule} belongs to
$\mathcal K_{i,+}$, induces $\Phi_\rho$, and globally maximizes
$\mathcal J_\rho$ over $\mathcal K_{i,+}$. Moreover,
\[
 0\le k_\rho(y)<K,\qquad y>0.
\]
All maximizing schedules induce the same investment policy from $(x,c)$.
More precisely, every maximizer $\tilde k$ satisfies $\tilde k(y)=k_\rho(y)$ and induces the response $\Phi_\rho(y)$ for $y>x\vee y^{\rm entry}_\rho$.
For $y<x\vee y^{\rm entry}_\rho$, no investment occurs from the
initial state, and subsidies need not be uniquely determined.
\end{proposition}

Appendix~\ref{app:regulatorExistenceProof} proves admissibility and
optimality, including comparison with capacity boundaries $\Phi(y)$ for which $M(y,\Phi(y))>K$
at some states. Thus the cap in \eqref{reg:candidateBoundary} is part
of the optimal construction, not an extra restriction on competitors.

 A finite subsidy cutoff exists
exactly when $\rho>(h\chi-1)/(1-\eta)$. At the cutoff demand $y_\rho^{\rm cut}$,
the relaxed and no-subsidy boundaries meet at capacity
\[
z_\rho^{\rm cut}
=\Phi_K(y_\rho^{\rm cut})
=\Phi_\rho^{\rm rel}(y_\rho^{\rm cut}).
\]
To see this condition, after substituting into \eqref{reg:regulatorFOC} $M(y,z)=K$  and re-arranging we obtain
\begin{equation}\label{reg:cutoffCapacity}
(\rho-\frac{h\chi-1}{1-\eta})z-\rho c
-\frac{\varepsilon}{rK(1-\eta)}G_A(z)=0.
\end{equation}
This capacity is the unique root $z>c$ and the corresponding cutoff is
\begin{equation}\label{reg:subsidyCutoff}
y_\rho^{\rm cut}
=\Phi_K^{-1}(z_\rho^{\rm cut})
=H^{-1}\!\left[
\frac{\chi K}{\eta}({z_\rho^{\rm cut})}^{1-\eta}\right].
\end{equation}
For $y\ge y_\rho^{\rm cut}$, $\Phi_\rho(y)=\Phi_K(y)$ and
$k_\rho(y)=0$.
When $\varepsilon=0$ and $\rho>(h\chi-1)/(1-\eta)$,
\[z_\rho^{\rm cut}=\frac{\rho(1-\eta)c}{\rho(1-\eta)+1-h\chi}
\]
When $\rho\le (h\chi-1)/(1-\eta)$, subsidies remain positive.

Finally, restore market indices and write
$k_\rho:=(k_{1,\rho},\ldots,k_{N,\rho})$. Thus,
\[
 S(k_\rho)=\sum_{i=1}^N S_i(k_{i,\rho}),
 \qquad \rho\ge\zeta.
\]

\begin{theorem}\label{thm:regulatorExistence}
The map $\rho\mapsto S(k_\rho)$ is continuous and nonincreasing,
tends to zero as $\rho\to\infty$, and satisfies
\[
 \rho_1<\rho_2,\quad S(k_{\rho_1})>0
 \quad\Longrightarrow\quad
 S(k_{\rho_2})<S(k_{\rho_1}).
\]
For every $B>0$, there is a unique $\rho^*\ge\zeta$ such that
\begin{equation}\label{reg:verifiedBudgetConditions}
 S(k_{\rho^*})\le B,\qquad
 (\rho^*-\zeta)[B-S(k_{\rho^*})]=0.
\end{equation}
\end{theorem}

The theorem links each positive budget to a unique common penalty
$\rho^*$ and the corresponding optimal nonnegative subsidy schedules.
A tighter binding budget raises this penalty and reduces total
subsidy expenditure. Once the budget is large enough, the constraint
no longer binds, and further increases in funding do not change
the optimal investment policy.

Finally, we conclude with a characterization of the shape of the optimal schedules

\begin{proposition}\label{prop:capacityBenefitShape}
For every $\rho\ge0$, define
\begin{equation}\label{reg:entryGap}
 \Delta_\rho
 :=k_\rho(y_\rho^{\rm entry})-(K-d_\rho).
\end{equation}
If $\Delta_\rho\le0$, the schedule $k_\rho$ is nonincreasing.
If $\Delta_\rho>0$, it is strictly increasing up to a unique
maximum $y_\rho^*>y_\rho^{\rm entry}$ and nonincreasing thereafter.

A sufficient condition for the first case is
\begin{equation}\label{reg:decreasingCondition}
 \frac{\varepsilon}{r}G_A'(c)(h\chi+\rho)
 \le \rho(1-\eta)(1+\rho)K.
\end{equation}
If $\varepsilon>0$, the second case holds for all sufficiently
small $\rho\ge0$.
\end{proposition}

\subsection{Zero capacity benefits and large budgets}
\label{subsec:constantSubsidies}

Suppose $\varepsilon_i=0$ in every market and the excess cost of raising funds is $\zeta=0$. Then the general
construction at $\rho=0$ gives $m_0=K/(h\chi)$ everywhere.
Consequently, the optimal subsidy is given by
\begin{equation}\label{reg:constantSubsidy}
 \ k^*_i=K_i\frac{h_i\chi_i-1}{h_i\chi_i}
\end{equation}
The required budget, including initial installation, is $B_{FB}=\sum_{i=1}^N k^*_i I_i(k^*_i)$

\section{Numerical Illustrations}
\label{sec:illustrations}

Throughout this section, we use $x$ as the demand variable and base the numerical simulations on a GBM demand process.
\[
dX(t)=\bar\mu X(t)dt+\bar\sigma X(t)dW(t).
\]
In that case, $\psi(x)=x^n$, where $n$ is the positive root of $\bar\sigma^2/2n(n-1)+\bar\mu n-r=0$.  Assumption~\ref{ass:diffusion} requires $r>\bar\mu$ while
Assumption~\ref{ass:regulatorExistence} additionally requires $n(1-\eta)>1$. This implies $n>1$ and hence that
$r>\bar\mu$.

\begin{remark}
\label{app:diffusionExamples}
Although the numerical examples use GBM demand, the results
also apply to more general diffusions. Examples with $X_0>0$ include:
\[
\begin{aligned}
\text{CIR:}\qquad
 \mathrm dX_t
 &=a(\theta-X_t)\dd t+\bar{\sigma}\sqrt{X_t}\dd W_t,\\
\text{Bessel:}\qquad
 \mathrm dX_t
 &=\frac{d_{\rm B}-1}{2X_t}\dd t+\dd W_t,\\
\text{Logistic:}\qquad
 \mathrm dX_t
 &=(aX_t-bX_t^2)\dd t+\bar{\sigma}X_t\dd W_t.
\end{aligned}
\]
Here $\bar\sigma,a,\theta,b>0$, with $2a\theta\ge \bar{\sigma}^2$
for CIR and $d_{\rm B}\ge2$ for Bessel demand.
For $r>0$ and $0<\eta<1$,
Assumption~\ref{ass:regulatorExistence} holds for CIR, Bessel,
and logistic demand. 
The special-function asymptotics used to verify
\eqref{reg:integrability} are given in
\citet[Chapters~10 and~13]{OlverEtAl2010}.
\end{remark}

We numerically compare optimal subsidies across heterogeneous markets,
budget levels, and capacity-benefit weights. The baseline
uses $r=0.12$, $\vartheta=1$, $\zeta=0$, $x_0=c_0=1$,
$(\alpha,\beta,\delta,K)=(0.45,0.5,3,1)$, and
$(\bar \mu,\bar \sigma)=(0.04,0.2)$. Initially, $\varepsilon_i=0$ in every
market. Budgets are measured relative to $B_{\mathrm{FB}}$, the expected
discounted expenditure needed to implement the unconstrained
welfare-maximizing investment policy.

\begin{figure}[!ht]
\centering
\includegraphics[
  width=0.93\linewidth,
  height=0.43\textheight,
  keepaspectratio
]{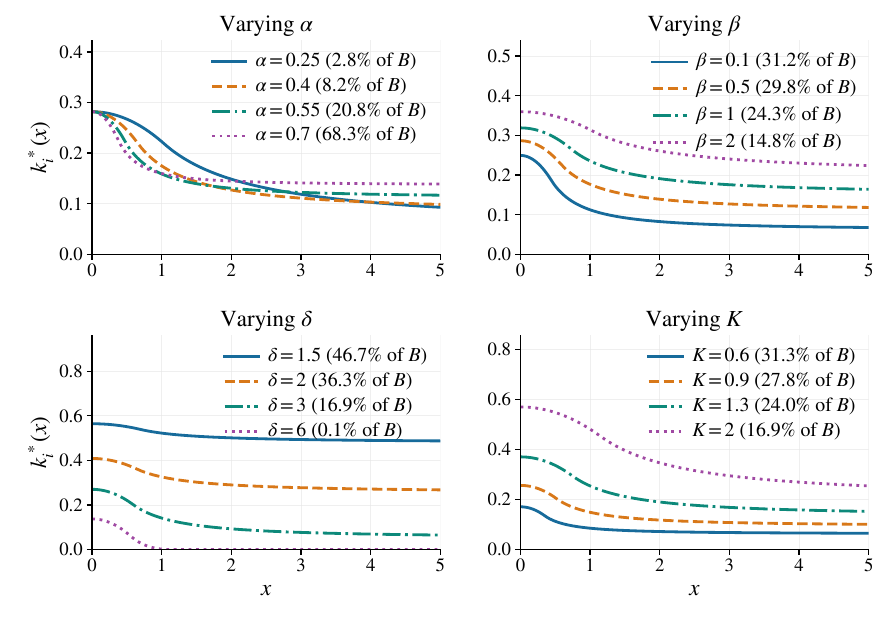}
\caption{Optimal subsidies across heterogeneous markets with
$\varepsilon_i=0$. Each panel varies one parameter across four markets
sharing a budget $B=0.2B_{\mathrm{FB}}$, calculated separately for that
economy. Legend percentages report expected discounted expenditure shares.}
\label{fig:gbm_comparative_statics}
\end{figure}

Figure~\ref{fig:gbm_comparative_statics} illustrates the importance of
tailoring support to market characteristics. In these examples, greater
demand elasticity is associated with lower subsidies and expenditure
shares, whereas higher installation costs or $\beta$ lead to larger
per-unit subsidies but smaller expenditure shares. Subsidy generosity
therefore does not directly measure a market's fiscal allocation.
Effective targeting must account for the timing and scale of induced
investment and the welfare gains it generates.

\begin{figure}[!htbp]
\centering
\includegraphics[
  width=0.5\linewidth,
  height=0.43\textheight,
  keepaspectratio
]{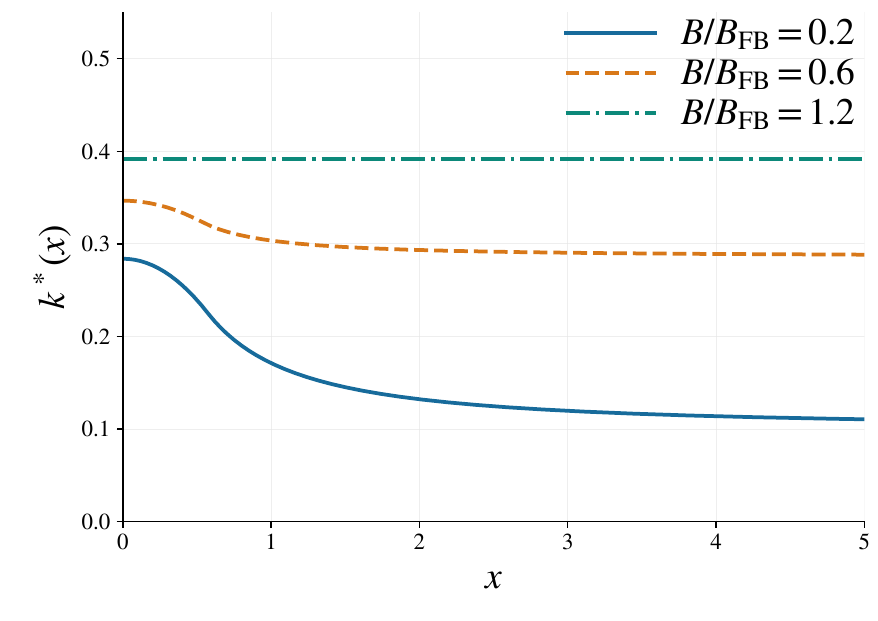}
\caption{Optimal subsidies across budget levels 
with $\varepsilon=0$ and $\zeta=0$.}
\label{fig:process_comparison}
\end{figure}

Figure~\ref{fig:process_comparison} shows how budget scarcity changes
policy design. With sufficient funding, the constant subsidy in
\eqref{reg:constantSubsidy} implements the unconstrained optimum.
Additional budget beyond $B_{\mathrm{FB}}$ leaves investment and
expenditure unchanged. When funding is scarce, optimal support is lower
and declines with demand.

\subsection{Positive capacity benefits}
\label{subsec:positiveBenefitNumerics}

\begin{figure}[!htbp]
\centering
\includegraphics[width=0.78\linewidth]{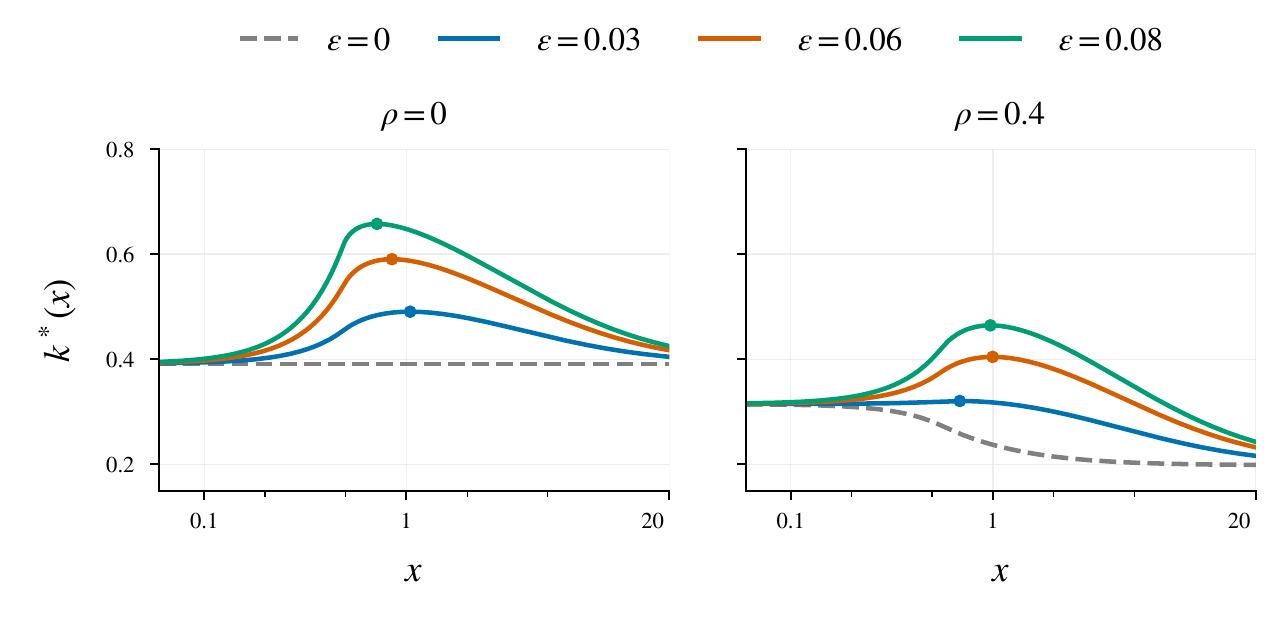}
\caption{Optimal subsidies with capacity benefits at $\rho=0$ (left)
and $\rho=0.4$ (right). Demand is shown on a logarithmic scale.
Dots mark subsidy maxima.}
\label{fig:positiveCapacityBenefits}
\end{figure}

We next use the capacity-benefit function
\begin{equation}
A(z)=\frac{b}{b-c}
\left[
b\log\left(\frac{b+z-c}{b}\right)
-c\log\left(\frac{z}{c}\right)
\right],
\qquad z\ge c,
\label{reg:rationalBenefit}
\end{equation}
which satisfies $A(c)=0$ and
$zA'(z)=b(z-c)/(b+z-c)$.
We set $b=25$, $c=K=1$, and
$\varepsilon\in\{0,0.03,0.06,0.08\}$, satisfying
Assumption~\ref{ass:capacityBenefit}.
Figure~\ref{fig:positiveCapacityBenefits} compares these weights
at fixed fiscal penalties however the implied budgets may differ.

Positive capacity benefits can change the direction of optimal support:
each positive-weight example has a subsidy that initially rises with
demand and then falls, consistent with
Proposition~\ref{prop:capacityBenefitShape}. Valuing capacity
for purposes such as strategic preparedness can therefore justify
increasing support over an initial range of demand, even when public
funds are scarce.

\appendix

\section{Proofs of Section~\ref{sec:producer}}

\begin{proof}[Proof of Corollary~\ref{cor:resolventProperties}]
Since $U$ in \eqref{eq:resolvent} is the expected discounted
integral of $X(t)$, \citet[Corollary~4.5, equation~(4.40)]
{LambertonZervos2013}, applied with running payoff $h(y)=y$,
gives
\[
\lim_{x\to\infty}\frac{U(x)}{\psi(x)}=0,
\]
because the upper endpoint is inaccessible.

To establish the properties of $H$, introduce the scale
and speed densities
\[
s'(x):=\exp\!\left(
-\int_{x_0}^x\frac{2\mu(y)}{\sigma^2(y)}\dd y
\right),
\qquad
m'(x):=\frac{2}{\sigma^2(x)s'(x)},
\]
where $x_0>0$ is arbitrary. Let $\varphi$ denote the positive
decreasing solution of $(\mathcal L-r)\varphi=0$ and note that the Wronskian divided by the scale density
\[
\mathcal W
:=\frac{\varphi(x)\psi'(x)
        -\varphi'(x)\psi(x)}{s'(x)}
\]
is a positive constant.

Applying their resolvent representation~(4.34) to the running
payoffs $h(y)=y$ and $h\equiv1$ gives
\[
\begin{aligned}
U(x)
&=\frac{\varphi(x)}{\mathcal W}
  \int_0^x y\psi(y)m'(y)\dd y
 +\frac{\psi(x)}{\mathcal W}
  \int_x^\infty y\varphi(y)m'(y)\dd y,\\
\frac1r
&=\frac{\varphi(x)}{\mathcal W}
  \int_0^x \psi(y)m'(y)\dd y
 +\frac{\psi(x)}{\mathcal W}
  \int_x^\infty \varphi(y)m'(y)\dd y.
\end{aligned}
\]

When differentiating either expression, the contributions
from the integration limits cancel. In particular,
\[
U'(x)
=\frac{\varphi'(x)}{\mathcal W}
  \int_0^x y\psi(y)m'(y)\dd y
 +\frac{\psi'(x)}{\mathcal W}
  \int_x^\infty y\varphi(y)m'(y)\dd y.
\]
Subtracting $q(x)$ times this derivative from $U(x)$
eliminates the integral over $(x,\infty)$, since
$q\psi'=\psi$. Moreover,
\[
\frac{\varphi(x)-q(x)\varphi'(x)}{\mathcal W}
=\frac{s'(x)}{\psi'(x)}.
\]
Thus, we obtain that
\[
H(x)
=\frac{s'(x)}{\psi'(x)}
  \int_0^x y\psi(y)m'(y)\dd y.
\]
The same calculation for $1/r$ gives
\[
\frac1r
=\frac{s'(x)}{\psi'(x)}
  \int_0^x \psi(y)m'(y)\dd y.
\]
Therefore
\[
H(x)
=\frac{\int_0^x y\psi(y)m'(y)\dd y}
       {r\int_0^x \psi(y)m'(y)\dd y}.
\]
Consequently, $0<H(x)<x/r$, and hence $H(0+)=0$.

Finally, the resolvent equation~(4.37) in the same reference
is $(\mathcal L-r)U=-x$.
Using this equation together with \eqref{eq:harmonic}
and the definitions in \eqref{eq:qH}, we obtain
\[
\begin{aligned}
H'(x)
&=q(x)\left(\frac{\psi''(x)}{\psi'(x)}U'(x)-U''(x)\right)\\
&=\frac{2q(x)}{\sigma^2(x)}(x-rH(x))>0.
\end{aligned}
\]
Direct differentiation also gives
\[
\left(\frac{U(x)}{\psi(x)}\right)'
=\frac{\psi(x) U'(x)-\psi'(x)U(x)}{\psi^2(x)}
=-\frac{\psi'(x)}{\psi^2(x)}H(x)<0.
\]
Finally, since $\lim_{x\rightarrow\infty}U(x)=\infty$, if $\lim_{x\rightarrow\infty}H(x)=M<\infty$, then
\[
U(x)=\psi(x)\int_x^\infty\frac{\psi'(y)}{\psi^2(y)}H(y)\,\dd y
\le M\psi(x)\int_x^\infty\frac{\psi'(y)}{\psi^2(y)}\,\dd y\le M,
\]
contradicting Assumption~\ref{ass:diffusion}.
Thus $H(x)\to\infty$ as $x\to\infty$.

\end{proof}

\begin{proof}[Proof of Proposition~\ref{prop:boundary}]
For each $c>0$, $\mathcal T(x,c)$ decreases continuously
from infinity to zero as $x$ increases. Hence there is exactly
one value $x=b(c)$ satisfying \eqref{eq:boundaryT}.
Increasing $c$ raises $\mathcal T(x,c)$ at every fixed $x$.
To restore the equality $\mathcal T(x,c)=\eta/\chi$,
$x$ must therefore increase implying that $b$ is strictly increasing.

The continuity, strict monotonicity, and endpoint limits
of $\mathcal T$ imply that
$\mathcal T(b(c),c)=\eta/\chi$ defines a unique
continuous, strictly increasing function $b$.
The capacity endpoint limits give $b(0+)=0$ and
$b(c)\to\infty$ as $c\to\infty$.

Direct differentiation  gives
\begin{equation}\label{eq:normalizedGainDerivative}
 \partial_x\!\left(\frac{g(x,c)}{\psi(x)}\right)
 =\frac{\psi'(x)}{\psi^2(x)}c^{\eta-1}H(x)
 \left[\mathcal T(x,c)-\frac{\eta}{\chi}\right].
\end{equation}
It is positive below $b(c)$ and negative above it. The boundedness of
$\ell$ and Corollary~\ref{cor:resolventProperties} imply that $g(x,c)>0$ for
large $x$ and $\lim_{x\rightarrow\infty}g(x,c)/\psi(x)=0$ proving that $g(x,c)>0$ on $x\geq b(c)$.
\end{proof}

\begin{proof}[Proof of Lemma~\ref{lem:HJB}]
We begin with the waiting region $\mathcal W$.
For $x<b(c)$, the candidate in \eqref{eq:candidate} is
\[
 w(x,c)=\frac{c^\eta}{\chi}U(x)+F(c)\psi(x).
\]
Since $U$ and $\psi$ are $C^2$ and satisfy
$(\mathcal L-r)U=-x$ and $(\mathcal L-r)\psi=0$,
\begin{equation}\label{eq:HJBWaiting}
 (\mathcal L-r)w(x,c)+\frac{xc^\eta}{\chi}
 =0.
\end{equation}
For the gradient constraint, differentiating \eqref{eq:F}
and using the definition of $g$ in \eqref{eq:marginalGain} gives
\begin{equation}
 w_c(x,c)-\ell(x,c)
 =\psi(x)\left[
   \frac{g(x,c)}{\psi(x)}
   -\frac{g(b(c),c)}{\psi(b(c))}
   \right]<0.
\end{equation}
The strict inequality follows from Proposition~\ref{prop:boundary}. Therefore, the HJB \eqref{eq:HJB} holds classically in the waiting region.

\medskip
In the investment region $\mathcal I$, differentiating
\eqref{eq:candidateLayers} and using
\eqref{eq:ConeMarginalValue} gives
\begin{equation}\label{eq:HJBGradient}
 w_c(x,c)-\ell(x,c)
 =g(x,c)-\mathcal Z(x,c)=0,
 \qquad x\ge b(c).
\end{equation}

Next, fix $c>0$ and put
\begin{equation}\label{eq:M}
 L(x,c):=q(x)g_x(x,c)-g(x,c)
 =c^{\eta-1}H(x)
   \left[\mathcal T(x,c)-\frac{\eta}{\chi}\right].
\end{equation}
By the boundary equation  \eqref{eq:boundaryT},
the strict decrease in \eqref{eq:thresholdSigns} giving $\mathcal T(x,c)<\eta/\chi$ for $x>b(c)$ and the fact that  Corollary~\ref{cor:resolventProperties} states that $H$ is
positive and nondecreasing we obtain that  $L(\cdot,c)$ is nonpositive and nonincreasing
there, and consequently of locally bounded variation in $x\ge b(c)$.

The identity
\[
 g_x(x,u)=\frac{L(x,u)+g(x,u)}{q(x)}
\]
also shows that $g_x(\cdot,c)$ has locally bounded variation for $x\ge b(c)$ since $g(\cdot,c)$ and $q$ are $C^1$.

Using $q=\psi/\psi'$ and dividing
$(\mathcal L-r)\psi=0$ by $\psi'$ gives
\begin{equation}\label{eq:qIdentity}
 \frac{\sigma^2(x)}{2q(x)}(1-q'(x))
       +\mu(x)=rq(x).
\end{equation}
The product rule for bounded-variation functions gives
\[
 \mathrm dL(x,u)
 =q(x)\,\mathrm d g_x(x,u)
      +(q'(x)-1)g_x(x,u)\,\mathrm dx.
\]
Substituting this identity into the definition of
$\mathcal G_r$ and using \eqref{eq:qIdentity}, we obtain,
for $x>b(u)$,
\begin{equation}\label{eq:factorization}
 \mathcal G_r[g(\cdot,u)](\mathrm dx)
 =\frac{\sigma^2(x)}{2q(x)}\,\mathrm dL(x,u)
       +rL(x,u)\,\mathrm dx <0.
\end{equation}
Indeed, $\mathrm dL(\cdot,u)\le0$ because $L(\cdot,u)$
is nonincreasing, and $L(\cdot,u)<0$.

For $x<b(c)$, $\mathcal{Z}(x,c)$ is a multiple of $\psi$ and is
therefore differentiable. At $x=b(c)$,
\eqref{eq:normalizedGainDerivative} and
\eqref{eq:boundaryT} imply
\[
 g_x(b(c),c)
 =\frac{\psi'(b(c))}{\psi(b(c))}g(b(c),c).
\]
Thus the first derivatives of the two branches of $\mathcal{Z}(\cdot,c)$
match at $b(c)$. It follows that $\mathcal{Z}(\cdot,c)$ is $C^1$,  $\mathcal{Z}_x(x,c)$ has locally bounded variation, and that
its second distributional derivative has no atom
at $b(c)$. Combining the two regions and this boundary
condition gives
\begin{equation}\label{eq:marginalGeneratorSign}
 \mathcal G_r[\mathcal{Z}(x,c)](\mathrm dx)\le0\quad\text{for }x>b(c),
 \qquad
 \mathcal G_r[\mathcal{Z}(x,c)](\mathrm dx)=0\quad\text{for }x\leq b(c).
\end{equation}

Fix $0<a<z$ and $Y>c\vee\Phi(z)$. For every $x\in[a,z]$
and $u\ge Y$, we have $x<b(u)$, so
\[
w(x,c)=\frac{c^\eta}{\chi}U(x)
+\int_c^Y\mathcal Z(x,u)\dd u+F(Y)\psi(x).
\]
Since $\mathcal Z_x$ is continuous, differentiation gives that $w(\cdot,c)\in C^{1}$. Together with
\eqref{eq:HJBGradient}, this gives $w\in C^1(\mathcal{D})$.

On the investment branch, $g_x=(L+g)/q$ and the properties
established above bound the variation of $g_x$ by a constant
independent of $u\in[c,Y]$.
On the waiting branch, $\mathcal Z_x$ is a multiple of
$\psi'$, with a coefficient bounded for $u\in[c,Y]$.
The branches match at $b(u)$, so $\mathcal{Z}_x(\cdot,u)$ is also of bounded variation. Integrating this bound over $[c,Y]$, and noting that the
remaining terms are $C^2$ in $x$, shows that
$w_x(\cdot,c)$ has bounded variation on $[a,z]$.

Since $(\mathcal L-r)U=-x$, we conclude that
\begin{equation}\label{eq:HJBResidualLayers}
 \mathcal G_r[w(\cdot,c)](\mathrm dx)
 +\frac{xc^\eta}{\chi}\,\mathrm dx
 =\int_c^\infty\mathcal G_r[\mathcal Z(\cdot,u)](\mathrm dx)\dd u
 \le0.
\end{equation}
Hence, together with \eqref{eq:HJBGradient}, this proves the HJB
in the measure sense in $\mathcal I$.
Together with the waiting-region calculation,
this establishes \eqref{eq:HJB} throughout $\mathcal D$.

\end{proof}

\begin{proof}[Proof of Theorem~\ref{thm:verification}]
For each fixed $c>0$, Proposition~\ref{prop:boundary}
and \eqref{eq:marginalGeneratorSign} give
\[
 \mathcal Z(\cdot,c)\ge\max\{g(\cdot,c),0\},
 \qquad
 \mathcal G_r[\mathcal Z(\cdot,c)]\le0.
\]
Thus, by the one-dimensional It\^o--Tanaka formula,
$e^{-rt}\mathcal Z(X_t,c)$ is a nonnegative local supermartingale.
Localizing and applying
Fatou's lemma shows that it is a supermartingale.
For any stopping time $\tau$, applying the same argument for the process stopped at $\tau$, together with $g\le\mathcal Z$ gives
\begin{equation}\label{eq:helper}
 \mathbb E_x\!\left[
 e^{-r\tau}g(X_\tau,c)\mathbf 1_{\{\tau<\infty\}}
 \right]\le\mathcal Z(x,c)
\end{equation} 

We apply this inequality to the installation times
generated by an arbitrary admissible capacity process
$C$, with $C_{0-}=c$. For each level $z>c$, let
\[
 \tau_z:=\inf\{t\ge0:C(t)\ge z\}
\]
If $C$ jumps at time $t$, all levels
$z\in(C(t-),C(t)]$ are installed at that time.
Their combined discounted expenditure is therefore
\[
 e^{-rt}\int_{C(t-)}^{C(t)}\ell(X(t),z)\,\mathrm dz.
\]
Tonelli's theorem and the strong Markov property give
\[
\begin{aligned}
 \mathbb E_{x,c}\!\left[
 \int_0^\infty e^{-rt}\frac{X(t)C(t)^\eta}{\chi}\,\mathrm dt
 \right]
 &=
 \frac{c^\eta}{\chi}U(x)
 +\int_c^\infty
 \mathbb E_{x,c}\!\left[
 e^{-r\tau_z}\frac{\eta}{\chi}z^{\eta-1}
 U(X(\tau_z));\,\tau_z<\infty
 \right]\mathrm dz\\
 &\quad\le
 w(x,c)+
 \mathbb E_{x,c}\!\left[
 \int_{[0,\infty)}e^{-rt}\ell(X(t),C(t))\,\mathrm dC_t
 \right]
 <\infty,
\end{aligned}
\]
where we used
$g(x,z)=(\eta/\chi)z^{\eta-1}U(x)-\ell(x,z)$
\eqref{eq:candidateLayers} and \eqref{eq:helper} and the fact that $w<\infty$.

Thus, $V\le w$.

To prove the reverse inequality, consider the candidate capacity process
\[
C^*(t)=c\vee\Phi(\overline X(t)).
\]
This process is adapted, nondecreasing, and finite on bounded time
intervals. Apart from a possible initial jump from $c$ to $C^*(0)$,
it is continuous and reflects $(X,C^*)$ at the boundary $x=b(c)$
in the positive capacity direction.

Indeed, by construction, $C^*(t)-\Phi(X(t))$ is nonnegative,
and $C^*$ increases only when this difference is zero.
Since $b$ is strictly increasing,
\[
X(t)\le b(C^*(t)),
\qquad
C^*(t)-C^*(0)
=
\int_{(0,t]}
\mathbf{1}_{\{X(s)=b(C^*(s))\}}\,\mathrm dC^*(s),
\qquad t\ge0.
\]
Moreover,
\[
X(t)\le\overline X(t)\le b(C^*(t))
=\max\{b(c),\overline X(t)\}.
\]
Whenever $X(t)=b(C^*(t))$, these inequalities imply
\[
X(t)=\overline X(t)=b(C^*(t)),
\]
The integral identity above shows that all capacity increases
after time zero occur at such times.

Next, levels
$c<z\le \Phi(x)$ are installed immediately, while each level
$z>c\vee \Phi(x)$ is installed when $X$ first reaches $b(z)$. The hitting-time identity
\eqref{eq:harmonic} gives

\begin{equation} \label{eq:verificationCost}
\mathbb E_{x,c}\!\left[
 \int_{[0,\infty)}e^{-rt}\ell(X(t),C^*(t))\,\mathrm dC^*(t)
 \right]=
 \int_c^{c\vee\Phi(x)}\ell(x,z)\,\mathrm dz
 +\psi(x)\int_{c\vee\Phi(x)}^\infty
 \frac{\ell(b(z),z)}{\psi(b(z))}\,\mathrm dz
 <\infty.
\end{equation}
Finiteness follows from Assumption~\ref{ass:tail} and $
 0<\ell(b(z),z)
 <\eta z^{\eta-1}U(b(z))/\chi$,
so $C^*$ is admissible.

For each level $z>c$, let $\tau_z$ denote its installation
time under $C^*$. We show that
\[
 \mathbb E_x\!\left[
 e^{-r\tau_z}g(X(\tau_z),z);\,\tau_z<\infty
 \right]
 =\mathcal Z(x,z).
\]

Suppose that $\Phi(x)>c$. Then $b(\Phi(x))=x$.
For every initially installed level $c<z\le \Phi(x)$,
monotonicity of $b$ gives $b(z)\le x$. Since $\tau_z=0$,
the definition of $\mathcal Z$ therefore yields
\[
 \mathbb E_x\!\left[
 e^{-r\tau_z}g(X(\tau_z),z);\,\tau_z<\infty
 \right]
 =g(x,z)
 =\mathcal Z(x,z).
\]

For every remaining level $z>c\vee \Phi(x)$, we have $x<b(z)$ and $
 \tau_z=\inf\{t\ge0:X(t)\ge b(z)\}$.
 Consequently, \eqref{eq:harmonic} gives
\begin{equation}
 \mathbb E_x\!\left[
 e^{-r\tau_z}g(X(\tau_z),z);\,\tau_z<\infty
 \right]=
 \frac{\psi(x)}{\psi(b(z))}g(b(z),z)=\mathcal Z(x,z).
\end{equation}

Applying the strong Markov
argument used above, the expected payoff of $C^*$ is
\[
\begin{aligned}
 &\mathbb E_{x,c}\!\left[
 \int_0^\infty e^{-rt}
 \frac{X(t)(C^*(t))^\eta}{\chi}\,\mathrm dt
 -
 \int_{[0,\infty)}e^{-rt}
 \ell(X(t),C^*(t))\,\mathrm dC_t^*
 \right]\\
 &\qquad=
 \frac{c^\eta}{\chi}U(x)
 +\int_c^\infty
 \mathbb E_x\!\left[
 e^{-r\tau_z}g(X(\tau_z),z);\,\tau_z<\infty
 \right]\mathrm dz\\
 &\qquad=w(x,c),
\end{aligned}
\]
where the last equality follows from
\eqref{eq:candidateLayers}. This establishes
$V(x,c)\ge w(x,c)$. Combined with $V(x,c)\le w(x,c)$
gives $V=w$.

Finally, since any optimal control must attain equality in \eqref{eq:helper} for
almost every capacity level $z>c$, installing before
$X$ reaches $b(z)$ gives strictly less than the value
of waiting. Any positive delay after $X$ first reaches
or exceeds $b(z)$ makes it spend positive time above
$b(z)$, and therefore causes a strict
loss by \eqref{eq:factorization}.
Thus each such level must be installed at
$\inf\{t\ge0:X(t)\ge b(z)\}$. These installation times
determine $C(t)=c\vee\Phi(\overline X(t))$.
By right-continuity, this equality holds for all $t\geq 0$ almost surely proving uniqueness.

\end{proof}

\section{Proofs of Section~\ref{sec:regulator}}
\label{app:regulatorExistenceProof}

\begin{lemma}\label{lem:regulatorCandidateAdmissibility}
The candidate
$\Phi_\rho:(0,\infty)\to(0,\infty)$ is a continuous, strictly increasing
bijection. Formula~\eqref{reg:candidateOptimalSchedule} defines an
admissible schedule $k_\rho\in\mathcal K_{i,+}$ that induces $\Phi_\rho$.
\end{lemma}
\begin{proof}
The derivative of the left-hand side of \eqref{reg:scalarThresholdMap} is
\[
 z^{-\eta}\left[(1-\eta)(1+\rho)K
 -\frac{\varepsilon}{r}
 \left(G_A'(z)-\eta\frac{G_A(z)}{z}\right)\right]
 \ge z^{-\eta}\left[(1-\eta)K
 -\frac{\varepsilon}{r}G_A'(c)\right]>0.
\]
Since $A'(z)\to0$, this side tends to infinity. The right-hand side
is nonincreasing in $z$ and strictly increasing in $y$, with equality
at $z=c$ exactly when $y=y^{\text{entry}}_\rho$. Thus $\Phi_\rho^{\rm rel}(y)$ is
continuous and strictly increasing from $c$ to infinity.
By \eqref{reg:scalarObjectiveDerivative}, for $y>y^{\text{entry}}_\rho$,
$F_\rho(y,\cdot)$ increases up to $\Phi_\rho^{\rm rel}(y)$ and decreases
thereafter, establishing its unique maximum while for $y\le y^{\text{entry}}_\rho$,
it decreases on $(c,\infty)$.

Define
\begin{equation}\label{reg:lowerAdjustedCost}
 \underline m:=
 \frac{(1+\rho)K-\frac{\varepsilon}{r}G_A'(c)}
      {h\chi+\rho}>0.
\end{equation}
Using $0\le A'(z)\le G_A'(c)$ in \eqref{reg:regulatorFOC},
together with \eqref{reg:candidateBoundary}, gives
\[
 \underline m\le m_\rho(y)\le K,
 \qquad
 \left(\frac{\eta H(y)}{\chi K}\right)^{1/(1-\eta)}
 \le\Phi_\rho(y)\le
 \left(\frac{\eta H(y)}{\chi\underline m}\right)^{1/(1-\eta)}.
\]
Both functions in the maximum defining $\Phi_\rho$ are strictly
increasing. At $y=y^{\text{entry}}_\rho$, the lower branch equals $c$, while the
upper branch converges to $\max\{c,\Phi_K(y^{\text{entry}}_\rho)\}=c$.
Hence, it is continuous. Next, by Corollary~\ref{cor:resolventProperties},
$H(y)\to0$ as $y\downarrow0$ and $H(y)\to\infty$ as
$y\to\infty$. The bounds above therefore show that
$\Phi_\rho$ is a
continuous strictly increasing bijection.

Since $\psi(y)\int_y^\infty\omega(t)\dd t=1$,
\eqref{reg:candidateOptimalSchedule} gives $
 0\le k_\rho(y)\le K-\underline m$ and
continuity of $m_\rho$ permits differentiation which implies that $k_\rho\in C^1$.
The producer's threshold function is
$\mathcal T(y,u)=(\eta/\chi)[u/\Phi_\rho(y)]^{1-\eta}$, which satisfies the
monotonicity and endpoint conditions in Assumption~\ref{ass:subsidy}.

Next, for every $v>0$, Tonelli's theorem and the bound on $\Phi_\rho$ give
\[
 \begin{aligned}
 \int_v^\infty u^{\eta-1}
 \frac{U(b_\rho(u))}{\psi(b_\rho(u))}\dd u
 &=\frac1\eta\int_{b_\rho(v)}^\infty
 \omega(y)H(y)[\Phi^\eta_\rho(y)-v^\eta]\dd y\\
 &\le\frac1\eta
 \left(\frac{\eta}{\chi\underline m}\right)^{\eta/(1-\eta)}\int_{b_\rho(v)}^\infty\omega(y)H^{1/(1-\eta)}(y)\dd y<\infty,
 \end{aligned}
\]
where finiteness follows from \eqref{reg:integrability}.

Finally, similar calculations show
\[
\begin{aligned}
S(k_\rho)
&\le K\,\mathbb E_{x,c}\!\left[
\int_{[0,\infty)} e^{-rt}\,dC(t)
\right]
=K\psi(x)\int_x^\infty
\omega(y)[\Phi_\rho(y)-c]^+\,dy
<\infty,
\end{aligned}
\]
Hence Assumption~\ref{ass:tail} holds, and
$k_\rho\in\mathcal K_{i,+}$ induces $\Phi_\rho$.
\end{proof}

The cap requires a global argument because $k\ge0$ does not imply
$M(y,\Phi(y))\le K$ pointwise. The following multiplier handles that distinction. 

\begin{lemma}\label{lem:regulatorPointwiseComparison}
Fix  $\rho\ge0$ and put
\[
 R_\rho(z):=(\rho-\frac{h\chi-1}{1-\eta})z-\rho c-\frac{\varepsilon}{rK(1-\eta)}G_A(z),\qquad z\ge c.
\]
Define
\begin{equation}\label{app:stateMultiplier}
 \lambda_\rho(y)=
 [R_\rho(c\vee\Phi_K(y))]_+.
\end{equation}
This function is nonnegative, nondecreasing, locally absolutely
continuous, and zero near the origin. It satisfies
$\lambda_\rho(y)[m_\rho(y)-K]=0$. For every $y,z>0$,
\begin{equation}\label{app:pointwiseComparison}
 F_\rho(y,z)-F_\rho(y,\Phi_\rho(y))
 \le\lambda_\rho(y)[M(y,z)-K],
\end{equation}
with strict inequality whenever $y>y^{\text{entry}}_\rho$ and $z\ne\Phi_\rho(y)$.
\end{lemma}
\begin{proof}
At $z=\Phi_K(y)$, using $M(y,z)=K$ and $A'(z)=G_A(z)/z$, we obtain
\[
\partial_zF_\rho(y,z)
=-\frac{K(1-\eta)}{z}
R_\rho(z).
\]
$R_\rho$ is convex, with $R_\rho(c)=-(h\chi-1)/(1-\eta)c<0$.
If $\rho\le \frac{h\chi-1}{1-\eta}$, it is negative everywhere. If $\rho>(h\chi-1)/(1-\eta)$, then
$\lim_{z\rightarrow\infty}R_\rho(z)=\infty$, so it has a unique zero and is strictly
increasing above that zero. Thus \eqref{app:stateMultiplier} defines
a nonnegative, nondecreasing, locally absolutely continuous function
$\lambda_\rho$ that is zero near the origin.

Above entry, Lemma~\ref{lem:regulatorCandidateAdmissibility}
and the preceding derivative identity imply
$\lambda_\rho(y)>0$ exactly when
$\Phi_\rho^{\rm rel}(y)<\Phi_K(y)$, in which case
$m_\rho(y)=K$. Below entry, $\lambda_\rho(y)=0$.
Hence
\[
\lambda_\rho(y)[m_\rho(y)-K]=0.
\]

Fix $y$ and set
$G_y(z):=F_\rho(y,z)-\lambda_\rho(y)M(y,z)$. For $z>c$,
\[
 G_y'(z)=M(y,z)\left[h\chi+\rho\eta+
 \frac{(1-\eta)(\rho c+\lambda_\rho(y))}{z}\right]
 +\frac{\varepsilon}{r}A'(z)-(1+\rho)K.
\]
As in \eqref{reg:scalarThresholdMap}, $G_y'$ has at most one zero,
with positive sign before it and negative sign afterwards.

Suppose first that $\lambda_\rho(y)>0$. On $z>c$,
$\Phi_\rho(y)=\Phi_K(y)>c$, and direct substitution gives
$G_y'(\Phi_\rho(y))=0$ and hence $\Phi_\rho(y)$ is
its unique maximizer. On the other hand for $0<z\le c$, we have
$G_y(z)=-\lambda_\rho(y)M(y,z)$, which is strictly increasing
because $M_z(y,z)<0$. Since $G_y$ is continuous at $c$, it is therefore strictly increasing on $(0,\Phi_\rho(y))$.
Thus $\Phi_\rho(y)$ is the unique global maximizer.

If instead $\lambda_\rho(y)=0$, then $G_y=F_\rho$.
The preceding characterization of $F_\rho$ shows that
$\Phi_\rho(y)$ is a global maximizer, and that this maximizer
is unique whenever $y>y^{\text{entry}}_\rho$.

Consequently, $G_y(z)\le G_y(\Phi_\rho(y))$, or equivalently,
\[
 F_\rho(y,z)-F_\rho(y,\Phi_\rho(y))
 \le \lambda_\rho(y)[M(y,z)-m_\rho(y)]
 =\lambda_\rho(y)[M(y,z)-K],
\]
where the last equality follows from
$\lambda_\rho(y)[m_\rho(y)-K]=0$.
This proves \eqref{app:pointwiseComparison}, with strict inequality
when $y>y^{\text{entry}}_\rho$ and $z\ne\Phi_\rho(y)$.
\end{proof}

\begin{proof}[Proof of Proposition~\ref{prop:regulatorCandidateVerification}]
Admissibility and the subsidy bounds follow from
Lemma~\ref{lem:regulatorCandidateAdmissibility}. Since $H$ is increasing,
\eqref{reg:revenueTailIdentity} and \eqref{reg:integrability} give
\begin{equation}\label{app:tailBoundaryVanishes}
 \frac{H^{1/(1-\eta)}(R)}{\psi(R)}
 \le\int_R^\infty\omega(y)H^{1/(1-\eta)}(y)\dd y\longrightarrow0.
\end{equation}

Note that $R_\rho(Z)\leq(\rho-(h\chi-1)/(1-\eta))Z$ and thus by \eqref{reg:candidateBoundary}
we have that $\lambda_\rho(y)=O(H^{1/(1-\eta)}(y))$. Consequently this implies that $\lim_{y\rightarrow\infty}\lambda_\rho(y)/\psi(y)=0$ and
$\int_x^\infty\lambda_\rho(y)\omega(y)dy<\infty$.
Integration by parts also gives $\int_x^\infty\lambda_\rho'(y)/\psi(y)dy<\infty$.

For any $k\in\mathcal K_{i,+}$ inducing the capacity boundary $\Phi(y)$
\[
 \omega(y)[M(y,\Phi(y))-K]=\left(\frac{k(y)}{\psi(y)}\right)'.
\]

Substituting $z=\Phi(y)$ into \eqref{app:pointwiseComparison},
multiplying by $\psi(x)\omega(y)$, and integrating over $[x,R]$,
we integrate by parts using the preceding identity.
Letting $R\to\infty$ and recalling
\eqref{reg:boundaryObjectiveRepresentation}, we obtain
\[
 \mathcal J_\rho(k)-\mathcal J_\rho(k_\rho)
 \le-\lambda_\rho(x)k(x)
 -\psi(x)\int_x^\infty
       \frac{\lambda_\rho'(y)k(y)}{\psi(y)}\dd y\le0.
\]
Boundedness of $k$ and \eqref{app:tailBoundaryVanishes} remove the upper
endpoint term. The signs follow from $k,\lambda_\rho,\lambda_\rho'\ge0$.

Now consider another maximizer $\tilde{k}$. Since it gives the same
objective value as $k_\rho$, equality must hold in
\eqref{app:pointwiseComparison} almost everywhere.
That inequality is strict when the boundaries differ above entry.
Continuity therefore gives $\Phi(y)=\Phi_\rho(y)$ for
$y>x\vee y^{\text{entry}}_\rho$, and
\eqref{reg:boundaryScheduleRecovery} gives $\tilde{k}(y)=k_\rho(y)$ there.
If $x<y^{\text{entry}}_\rho$, both boundaries remain below $c$
before entry, so neither induces investment.
Thus both schedules give the same investment policy from $(x,c)$.
\end{proof}

\begin{proof}[Proof of Theorem~\ref{thm:regulatorExistence}]
Fix a market. By \eqref{reg:lowerAdjustedCost},
\[
 m_\rho(y)\ge\underline m\ge
 \frac{K-\frac{\varepsilon}{r}G_A'(c)}{h\chi}>0,
 \qquad
 \Phi_\rho(y)\le
 \left(\frac{\eta hH(y)}
 {K-\frac{\varepsilon}{r}G_A'(c)}\right)^{1/(1-\eta)}.
\]
These bounds do not depend on $\rho$.
Equations~\eqref{reg:scalarThresholdMap} and
\eqref{reg:candidateBoundary} give continuity of $\Phi_\rho(y)$
and $m_\rho(y)$ in $\rho$. The bounds above and
\eqref{reg:integrability} therefore justify passing limits through
the expenditure integral \eqref{reg:boundaryTransferRepresentation}.
Hence $S(k_\rho)$ is continuous. The same capacity bound in
\eqref{reg:stateInvestment} bounds $I(k_\rho)$ uniformly.

To show that spending vanishes, the subsidy formula
\eqref{reg:candidateOptimalSchedule} gives
\[
 0\le k_\rho(y)\le K-\underline m
 =\frac{K(h\chi-1)+\frac{\varepsilon}{r}G_A'(c)}
        {h\chi+\rho}.
\]
Consequently,
$S(k_\rho)\le(K-\underline m)I(k_\rho)\to0$.

For $\rho_1<\rho_2$, each schedule maximizes the objective at its
own penalty, so
\[
\mathcal J_{\rho_1}(k_{\rho_1})
\ge \mathcal J_{\rho_1}(k_{\rho_2}),
\qquad
\mathcal J_{\rho_2}(k_{\rho_2})
\ge \mathcal J_{\rho_2}(k_{\rho_1}).
\]
Since $\rho$ enters the objective only through $-\rho S(k)$,
adding these inequalities gives
\[
(\rho_2-\rho_1)
[S(k_{\rho_1})-S(k_{\rho_2})]\ge0.
\]
Thus spending cannot increase. If it remained equal and positive,
both schedules would maximize $\mathcal J_{\rho_1}$.
By Proposition~\ref{prop:regulatorCandidateVerification}, their
boundaries would agree wherever investment occurs. Positive
spending then gives an interval with a common capacity $z>c$
and adjusted cost $m<K$. Subtracting
\eqref{reg:regulatorFOC} at the two penalties would require
\[
 K=m\left[\eta+(1-\eta)\frac cz\right]<m<K,
\]
a contradiction since $\eta<1$. Spending therefore decreases strictly while positive.

These conclusions also hold after summing across markets.
If $S(k_\zeta)\le B$, set $\rho^*=\zeta$.
Otherwise, continuity, strict decrease and the zero limit give
a unique $\rho^*>\zeta$ with $S(k_{\rho^*})=B$.
This proves \eqref{reg:verifiedBudgetConditions}, and
Proposition~\ref{prop:regulatorVerification} establishes optimality.

Finally, let $\widetilde k$ be another feasible optimum and put
$\Lambda^*=\rho^*-\zeta$. Equal objective values and
\eqref{reg:verifiedBudgetConditions} imply
\[
 0=\sum_i\left[
 \mathcal J_{i,\rho^*}(k_i^*)-
 \mathcal J_{i,\rho^*}(\widetilde k_i)\right]
 +\Lambda^*[B-S(\widetilde k)].
\]
Each term is nonnegative, so each must be zero.
Thus every $\widetilde k_i$ also maximizes $\mathcal J_{i,\rho^*}$,
and Proposition~\ref{prop:regulatorCandidateVerification}
gives the stated uniqueness.
When $B = 0$, every feasible response receives zero expected subsidy, so its producer payoff equals its payoff under zero subsidies.

\end{proof}

\begin{proof}[Proof of Proposition~\ref{prop:capacityBenefitShape}]
Recall that $m_\rho=d_\rho$ below entry. Differentiating \eqref{reg:regulatorFOC} and using
concavity of $G_A$ for $z=\Phi^{\rm rel}_\rho(y)$ shows that $m_\rho$ is either
nondecreasing or first decreasing and then nondecreasing.
The first case holds under \eqref{reg:decreasingCondition}.
In the second case, an increase must follow the initial
decrease, since
\[
 m_\rho(y_\rho^{\rm entry})=d_\rho,\qquad
 \lim_{y\to\infty}m_\rho(y)
 =K\wedge\frac{(1+\rho)K}{h\chi+\rho\eta}
 \ge d_\rho.
\]

Differentiating \eqref{reg:candidateOptimalSchedule}
and integrating by parts gives
\begin{equation}\label{app:shapeTail}
 \frac{k_\rho'(y)}{\psi'(y)}
 =\frac{k_\rho(y)-K+m_\rho(y)}{\psi(y)}
 =-\int_y^\infty\frac{m_\rho'(t)}{\psi(t)}\dd t,
\end{equation}
where $m'_\rho$ is an almost everywhere derivative since entry and cutoff joins may create kinks. The boundary term vanishes because $m_\rho$ is bounded
and $\psi\to\infty$.
Below entry, this ratio equals
$\Delta_\rho/\psi(y_\rho^{\rm entry})$.
While $m_\rho$ decreases, the ratio strictly decreases,
since its derivative is $m_\rho'/\psi$.
On the subsequent nondecreasing segment,
\eqref{app:shapeTail} makes it nonpositive.
Thus $\Delta_\rho\le0$ implies $k_\rho'\le0$.
The same conclusion holds whenever $m_\rho$ is
nondecreasing, proving \eqref{reg:decreasingCondition}.

If $\Delta_\rho>0$, $m_\rho$ must initially decrease.
The ratio in \eqref{app:shapeTail} then decreases
strictly from a positive value at entry to a negative
value at the end of that segment, because a subsequent
increase of $m_\rho$ remains in the tail integral.
It therefore crosses zero exactly once on that segment,
at some $y_\rho^*>y_\rho^{\rm entry}$, and is nonpositive
thereafter. Hence $k_\rho$ strictly increases up to $y_\rho^*$ and is nonincreasing
afterwards.

Finally, if $\varepsilon>0$, positivity of $G_A(z)$
for $z>c$ gives $A'(z)>0$, so
\eqref{reg:regulatorFOC} implies $m_0(y)<d_0$
above entry. Therefore,
\[
 \Delta_\rho
 =\psi(y_\rho^{\rm entry})
  \int_{y_\rho^{\rm entry}}^\infty
       [d_\rho-m_\rho(t)]\omega(t)\dd t
\]
yields $\Delta_0>0$.
Continuity and dominated convergence
give $\Delta_\rho\to\Delta_0$ as $\rho\downarrow0$.
Hence $\Delta_\rho>0$ for all sufficiently small $\rho\ge0$.
\end{proof}

\bibliographystyle{plainnat}
\bibliography{references}

\end{document}